\documentclass{ifacconf}

\usepackage{graphicx}      
\usepackage{natbib}        
\usepackage{amssymb}
\usepackage{xcolor}
\usepackage{amsmath}
\usepackage{makecell}
\newtheorem{lemma}{Lemma}
\newtheorem{assumption}{Assumption}
\newtheorem{remark}{Remark}

\usepackage{booktabs} 
\usepackage{multirow}

\begin{document}
\begin{frontmatter}

\title{A Communication-Efficient Distributed Optimization Algorithm for Problems with Coupling Constraints \thanksref{footnoteinfo}} 

\thanks[footnoteinfo]{This work was supported in part by National Key R\&D Program	of China under the grant 2022YFB3303900.}

\author[First]{Yuzhu Duan}
\author[First]{Ziwen Yang} 
\author[First]{Xiaoming Duan} 
\author[First]{Shanying Zhu}

\address[First]{School of Automation and Intelligent Sensing, Shanghai Jiao Tong University, Shanghai 200240, China.\\
 Key Laboratory of System Control and Information Processing, Ministry of Education of China, Shanghai 200240, China.\\
 Shanghai Key Laboratory of Perception and Control in Industrial Network Systems, Shanghai 200240, China.}

\begin{abstract}             
Resource allocation is a fundamental problem in Industrial Internet of Things (IIoT) systems, in which devices work together under limited communication bandwidth to complete diverse tasks. This paper proposes a communication-efficient distributed optimization algorithm tailored for problems with coupled constraints. To tackle coupled constraints, we solve the problem via  its dual counterpart, and develop a compressed version.  Difference compression and dynamic scaling factors are then introduced to mitigate compression errors. We show that the proposed algorithm converges linearly for strongly convex and smooth objective functions.  Numerical simulations validate the theoretical results and demonstrate the efficiency and robustness of the proposed  algorithm. 
	
\end{abstract}


\end{frontmatter}

\section{Introduction} \label{section1}

The Industrial Internet of Things (IIoT)  is becoming a new paradigm for industrial production environments.
Within industrial production settings, key issues including production resource planning \citep{erol2012multi}, task allocation for mobile intelligent agents \citep{nishi2005distributed}, and logistics management \citep{zhang2018framework} can be formulated as resource allocation problems,  where inherent characteristics of task-resource coupling within and across nodes introduce complex coupled constraints. 
Numerous distributed optimization methods have been developed for resource allocation problems in industrial settings \citep{guo2020adaptive, lee2021fast}. While these approaches enable agents to collaboratively minimize a global objective via local coordination \citep{ding2021differentially, xu2018bregman}, they face escalating communication pressures from the growing number of devices. 

To address the bottleneck caused by limited communication bandwidth, an effective approach is communication compression, with sparsification and quantization being common techniques \citep{karimireddy2019error, alistarh2017qsgd, zhu2018mitigating}. However, compression operators inevitably introduce quantization errors, which can degrade algorithm convergence performance—for instance, convergence only to a neighborhood of the optimal solution \citep{nedic2008distributed, aysal2008distributed}. To mitigate the impact of compression errors and improve the accuracy of convergence, researchers have proposed the compressed gradient difference scheme \citep{mishchenko2025distributed, tang2019doublesqueeze}. It is noted that these studies \citep{mishchenko2025distributed, tang2019doublesqueeze} were implemented under a master-worker framework, which limits their applicability to systems with a large number of devices.
On the other hand, variable difference compression scheme has been proposed to achieve efficient communication \citep{tang2018communication}. Novel algorithms with such scheme were proposed in \citep{koloskova2019decentralized1, koloskova2019decentralized}. These algorithms combined distributed gradient descent with model averaging. However, they achieved only a sublinear convergence rate, even for strongly convex objective functions.

To achieve faster convergence in the presence of compression, researchers have developed  several distributed optimization algorithms. Under relative compression errors, compressed optimization algorithms \citep{liao2022compressed, song2022compressed} incorporating gradient tracking technique have been proposed, yielding a linear convergence rate. The authors in  \citep{xiong2022quantized} proposed a distributed algorithm under absolute compression errors, achieving linear convergence. However, this approach, along with those in \citep{liao2022compressed, song2022compressed, Duan2024ASCC}, is ineffective for achieving communication efficiency under coupled constraints. 

To enhance communication efficiency for distributed optimization under coupled constraints, the authors in \citep{Qian2025} developed a dual gradient-tracking algorithm.
Its linear convergence was only guaranteed for uniform quantization and could not be extended to general quantizers. Meanwhile, in \citep{shi2024distributed}, a push–pull distributed algorithm under a uniform quantizer with finite quantization levels was proposed, which achieved linear convergence for strongly convex and smooth problems subject to local and global coupled constraints.  Nevertheless, these existing methods fail to characterize intricate coupling relationships, especially the internal coupling relations within individual nodes, nor could they offer convergence guarantees for a wider class of quantizers.

In this work, we aim to propose a communication-efficient distributed optimization algorithm capable of handling coupled constraints under various compressors, while guaranteeing linear convergence. The contributions  are summarized as follows. 
\begin{itemize}
\item[(1)] 
We propose a novel communication-efficient distributed optimization algorithm for constraint-coupled problems. 
To the best of our knowledge, few existing algorithms have achieved communication efficiency while handling coupled constraints under general compressors. See Table \ref{table} for comparison results.
\item[(2)] The proposed algorithm provably achieves linear convergence for minimizing strongly convex and smooth objective functions under both unbiased and biased compressors.  Numerical results demonstrate that the proposed compressed algorithm is robust under various compressors and remains effective in the presence of coupled equality constraints.
\end{itemize} 

\begin{table}[hb]
	\begin{center}
		\setlength{\tabcolsep}{1pt}
		\caption{Comparison with existing compressed distributed optimization algorithms.}\label{table}
		\begin{tabular}{cccccc}
			\hline \textbf{References} & \textbf{\makecell{Relative\\error}} & \textbf{\makecell{Absolute\\error}}  & \textbf{\makecell{Convergence\\rate}} & \textbf{\makecell{Coupled\\ Constraint}} \\\hline
			\makecell{Koloskova (2019a)\\Koloskova (2019b) } & $\checkmark$ & $\times$ & sublinear &  $\times$ \\
			\makecell{\citep{liao2022compressed}\\ \citep{song2022compressed}} &$\checkmark$ & $\times$  & linear &$\times$  \\
			\citep{xiong2022quantized} &$\times$  &$\checkmark$  & linear &$\times$  \\
			\citep{Duan2024ASCC} &$\checkmark$ & $\checkmark$ & linear &$\times$  \\
		\citep{ren2026distributed} &$\checkmark$& $\checkmark$& linear &$\times$  \\	\citep{Qian2025} & 	$\times$& $\checkmark$&linear &	$\checkmark$\\
		\citep{shi2024distributed} & 	$\times$& $\checkmark$&linear &	$\checkmark$\\
			\textbf{Our paper}   &$\checkmark$   &$\checkmark$    & linear  &$\checkmark$ \\ 
			\hline
		\end{tabular} 
	\end{center}
\end{table}

The remainder of this paper is organized as follows. Section \ref{section2} formulates the problem and introduces the compression model. The proposed compression algorithm is described in detail in Section \ref{section3}. Convergence analysis of the algorithm is provided in Section \ref{section4}. Simulation examples are presented in Section \ref{section5}. Finally, conclusions are drawn in Section \ref{section6}.

\textit{Notations:} Let $\mathbf{x}=\left[x_{1}^{T}, x_{2}^{T},..., x_{m}^{T}\right]^{T}$ denote the
collection of local variables $x_i$. We denote by $z_{i,k}$ and $\mathbf{z}_{k}$
the iterates of $z_i$ and $\mathbf{z}$ at time $k$. In addition, we use $\mathbf{1}$ to denote an  all-ones column vector. $\|\cdot\|$ represents the Euclidean norm of a vector, and $\Delta$ denotes the difference between two consecutive vectors, e.g., $\Delta \mathbf{z}_{k+1}=\mathbf{z}_{k+1}-\mathbf{z}_{k}.~\langle \cdot, \cdot \rangle$ is the inner product. $\otimes$ denotes the
Kronecker product.  $\mathbf{I}$ is the identity matrix with proper dimensions.  Let $\mathcal{H}$ denote the Euclidean space, and define a $\mathbf{G}$-space and its induced norm as $\langle \mathbf{z}, \mathbf{z}^{'}\rangle_{\mathbf{G}}=\langle\mathbf{G}\mathbf{z}, \mathbf{z}^{'}\rangle$ and 
$\|\mathbf{z}\|_{\mathbf{G}}=\sqrt{\langle \mathbf{G}\mathbf{z}, \mathbf{z}^{'}\rangle}, \forall~ \mathbf{z}, \mathbf{z}^{'} \in \mathcal{H}$, where $\mathbf{G}$ is a positive definite matrix. 
For some convex function $f(\cdot)$, its convex conjugate is denoted as $f^{*}(\mathbf{y}):=\sup _{ \mathbf{z}\in \mathcal{H}}\{\langle \mathbf{z}, \mathbf{y}\rangle-f(\mathbf{z})\}$.
 \(\mathbb{E}\left[\cdot\right]\) and \(\mathbb{E}\left[\cdot | \cdot\right]\) denote the expectation and conditional expectation, respectively.

\section{Problem Formulation and Preliminaries}
\label{section2}

\subsection{Distributed Optimization}
We consider a network with $m$ agents, where each agent has a local objective function $f_{i}:\mathbb{R}^{d} \rightarrow \mathbb{R}$. All agents solve the following optimization problem with coupled constraints:
\begin{equation} \label{problem}
	\begin{aligned} 
		&\min_{\mathbf{z} \in \mathbb{R}^{md}}&& f(\mathbf{z}) = \sum_{i = 1}^{m} f_{i}(z_{i})\\
		&\text{s.t.} &&\sum_{i=1}^{m} A_iz_{i}=\sum_{i=1}^{m} b_{i},
	\end{aligned}
\end{equation}
where  $\mathbf{z}=[z_1^{T},z_2^{T},...,z_{m}^{T}]^{T}\in \mathbb{R}^{md}$, $A_i\in \mathbb{R}^{n\times d} (n \le d)$ is the coupling matrix such that $\left[A_1,A_2,...,A_m\right]$ has full row
rank, and $\sum_{i=1}^{m} b_{i}\in \mathbb{R}^{n\times 1}$ is the load demand.  We make the following assumptions on the local objective functions:

\begin{assumption} \label{existence}
	There exists at least one  finite optimal solution to problem  \eqref{problem}.
\end{assumption}
\begin{assumption} \label{smooth}
	 $f_i:\mathbb{R}^{d}\rightarrow \mathbb{R}$ is $L_{f_i}$-Lipschitz smooth and $l_{f_i}$-strongly convex, i.e., for any $z_1, z_2\in \mathbb{R}^{d}$, 
\begin{subequations}
	\begin{equation}
	\|\nabla f_i(z_1)-\nabla f_i(z_2)\|_2\le L_{f_i}\|z_1-z_2\|_2,
		\end{equation}
	\begin{equation}
		(z_1-z_2)^{T}(\nabla f_i(z_1)-\nabla f_i(z_2))\ge l_{f_i} \|z_1-z_2\|_2^2,
		\end{equation}
\end{subequations}
where $L_{f_i}>0$ and $l_{f_i}>0$ are the Lipschitz and strong convexity constants, respectively.

It is not difficult to see that $f$ has $L_f$-Lipschitz gradient with $L_f=\max\{L_{f_i}\}$, and $f$ is $l_f$-strongly convex with $l_f=\min\{l_{f_i}\}$.
\end{assumption}

\begin{remark}
	Assumptions \ref{existence}-\ref{smooth} ensure the existence and uniqueness of the optimal solution $\mathbf{z}^{*}\in \mathbb{R}^{md}$ to  \eqref{problem}.
\end{remark}

\subsection{Basics of Graph Theory}

The exchange of information between agents is captured by an undirected graph $\mathcal{G}=(\mathcal{V}, \mathcal{E})$, where $\mathcal{V}=\{1,..,m\}$ is the set of agents and $\mathcal{E} \subseteq \mathcal{V}\times \mathcal{V}$ is the set of edges. $(i,j) \in \mathcal{E}$ if and only if agents $i$ and $j$ can communicate with each other. Let ${W}=[w_{ij}]\in \mathbb{R}^{m\times m}$ be the weight matrix of $\mathcal{G}$, namely $w_{ij}>0$ if $(i,j)\in \mathcal{E}$ or $i=j$, and $w_{ij}=0$ otherwise. Meanwhile, $\mathcal{N}_{i}=\{j\in \mathcal{V}|(i,j)\in \mathcal{E}\}$ denotes the neighbor set of agent $i$. 

\begin{assumption} \label{matrix}
The weight matrix ${W}$  satisfies the following conditions:
\begin{subequations}
	\begin{equation}
 \text{(Positive-definiteness)}~~{W}^{T}={W}~\text{and}~{W}\succ0,
	\end{equation}
	\begin{equation}
	 \text{(Stochasticity)}~~~ {W}\mathbf{1}=\mathbf{1} ~\text{or}~ \mathbf{1}^{T}{W}=\mathbf{1}^{T},
	\end{equation}
	\begin{equation}
		\text{(Connectivity)}~~~ \eta:=\rho\left({W}-\frac{\mathbf{1}\mathbf{1}^{T}}{m}\right)<1.	\end{equation}
\end{subequations}
\end{assumption}

\subsection{Compression Model}
 We introduce the following different types of compressors $Q(\cdot)$,
 \begin{itemize}
\item[(i)] For some $\sigma \in [0,1)$, the compressor $Q_1(\cdot)$ satisfies:
 \begin{equation}
 	\mathbb{E}\left[Q_1({x})\right]={x},~\mathbb{E}\left[\|Q_1({x})-x\|^2\right]\le\sigma^2, ~\forall {x} \in \mathbb{R}^{n}.
 \end{equation}	

 \item[(ii)] For some $C>0$, the compressor $Q_2(\cdot)$ satisfies:
  \begin{equation}
 	\mathbb{E}\left[Q_2(x)\right]=x,~\mathbb{E}\left[\|Q_2(x)-x\|^2\right]\le C\|x\|^2, ~\forall x\in \mathbb{R}^{n}.
  \end{equation}

 \item[(iii)] For some $\sigma \in [0,1)$,  the compressor $Q_3(\cdot)$ satisfies:
	\begin{equation}
		\|Q_3(x)-x\|_{2}^2\le \sigma^2,~x\in \mathbb{R}^{n}.
	\end{equation}
\end{itemize}

\begin{remark}
Compressors (i) and (ii)  are unbiased stochastic compression operators, which adopt absolute and relative compression errors for the input $x\in \mathbb{R}^{n}$, respectively.  Moreover,  the widely-used deterministic quantizers, such as compressor (iii), are biased compression operators with absolute compression error \citep{magnusson2020maintaining, xiong2022quantized}.
\end{remark}

The main objective of this paper is to design a distributed algorithm
where agents are only allowed to communicate compressed variables with their neighbors, with linear convergence to the exact optimal solution $\mathbf{z}^*\in \mathbb{R}^{md}$ of problem \eqref{problem} under different compressors.

\section{Communication-efficient Distributed Algorithm Design} \label{section3}

In this section, we design a communication-efficient distributed algorithm under different compressors and analyze its convergence properties.

To achieve a distributed solution, we first leverage duality theory to transform the original optimization problem \eqref{problem} into a  consensus problem as follows, 
\begin{equation}  \label{dualproblem}
	\begin{aligned} & \min_{\mathbf{x} \in \mathbb{R}^{nm}}&& \varphi(\mathbf{x})=\sum_{i=1}^{m} \varphi_{i}\left(x_{i}\right) \\ & s.t. && x_{i}=x_{j}, ~~ \forall i, j \in \mathcal{V}, \end{aligned}
\end{equation}
where $\varphi_i(x_i)=\sup_{z_i}\{-f_i(z_i)+x_i^T A_i z_i\}$ denotes the
local dual function.  Encoding $b_i$ in the initial value of $z_i$, the equivalent problem of \eqref{problem} can be written as
\begin{equation} \label{equi}
	\begin{aligned}
		&\min_{\mathbf{y} \in \mathbb{R}^{nm}}&& H(\mathbf{y})=\sum_{i=1}^{m} h_{i}\left(y_{i}\right),\\
		&s.t. &&\sum_{i=1}^{m} y_{i}=0,
	\end{aligned}
\end{equation}
where $h_i(y_i)=\inf_{z_i, A_i z_i=y_i} \{f_{i}(z_i)\}$. To tackle problem \eqref{problem}, 
\citep{wang2020dual} proposed the following distributed optimization algorithm:
\begin{subequations}  \label{DuSPA}
	\begin{equation} \label{DuSPA1}
		x_{i,k+1}=\sum_{j\in \mathcal{N}_{i}\cup\{i\}}w_{ij}x_{i,k}+\tau(y_{i,k}-A_i z_{i,k}),
	\end{equation}
	\begin{equation}  \label{DuSPA2}
		y_{i,k+1}=y_{i,k}-\frac{1}{\tau}\sum_{j\in \mathcal{N}_{i}}w_{ij} (x_{i,k+1}-x_{j,k+1}),
	\end{equation}
	\begin{equation}  \label{DuSPA3}
		z_{i,k+1}=z_{i,k}-\gamma \nabla f_{i}(z_{i,k})+\gamma A_{i}^{T}(2x_{i,k+1}-x_{i,k}),
	\end{equation}
\end{subequations}
where $x_{i,k}\in \mathbb{R}^{n}$ is the dual variable in \eqref{dualproblem}, $y_{i,k} \in \mathbb{R}^{n} $ is the auxiliary variable in \eqref{equi}, $\tau, \gamma$ are positive parameters.

To implement algorithm \eqref{DuSPA}, at each iteration, each agent $j$ needs to exactly communicate  $x_{i,k} \in \mathbb{R}^{n}$ and $x_{i,k+1} \in \mathbb{R}^{n}$ with its neighbors, which requires a significant amount of data exchange especially when the dimension $n$ is large. However, communication bandwidth is limited in practice. Here, we consider incorporating communication compression to resolve the problems caused by limited communication bandwidth.

\subsection{Algorithm Development}

To reduce communication overhead, each agent  transmits only compressed information $Q(\cdot)$ to its neighbors. Directly compression of the state variables leads to large compression errors. Inspired by DIANA \citep{mishchenko2025distributed} and LEAD \citep{liulinear}, we introduce an auxiliary variable $h_{i,k}$ as a reference point for  $x_{i,k}$ and compress the difference $x_{i,k}-h_{i,k}$ instead.  
However, for compressors (i) and (ii),  using only difference compression causes the absolute compression errors to dominate the signal as the difference decays to zero.  To address this issue, we introduce a dynamic scaling factor $r_k$  to rescale the difference into a proper range, preventing small signal from being overwhelmed by compression errors. Hence, leveraging the difference compression and dynamic scaling  technique,  the distributed optimization algorithm \eqref{DuSPA} can be modified to a communication-efficient version by incorporating the compressor
\begin{subequations} \label{compressalgorithm}
	\begin{equation}  \label{equation01}
		x_{i,k+1}=x_{i,k}-\psi (\hat{x}_{i,k}-\hat{x}_{W,i,k})+\tau (y_{i,k}-A_i z_{i,k}),
	\end{equation}
	\begin{equation}  \label{equation02}
		\tau y_{i,k+1}=\tau y_{i,k}-\psi (\hat{x}_{i,k+1}-\hat{x}_{W,i,k+1}),
	\end{equation}
	\begin{equation}  \label{equation03}
		z_{i,k+1}=z_{i,k}-\gamma \nabla f(z_{i,k}) +\gamma A_i^{T}(2 x_{i,k+1}-x_{i,k}),
	\end{equation}
	\text{with}
	\begin{equation}  \label{compress01}
		\hat{x}_{i,k}=h_{i,k}+r_{k}Q\left(\frac{x_{i,k}-h_{i,k}}{r_{k}}\right),
	\end{equation}
	\begin{equation}  \label{compress11}
		\hat{x}_{W,i,k}=\tilde{h}_{i,k}+\sum_{j\in \mathcal{N}_i}w_{ij}r_{k}Q\left(\frac{x_{j,k}-h_{j,k}}{r_{k}}\right),
	\end{equation}
	\begin{equation}  \label{compress02}
		h_{i,k+1}=(1-\alpha)h_{i,k}+\alpha \hat{x}_{i,k},
	\end{equation}
	\begin{equation}  \label{compress21}
		\tilde{h}_{i,k+1}=(1-\alpha)\tilde{h}_{i,k}+\alpha\hat{x}_{W,i,k},
	\end{equation}
\end{subequations}
where $\psi > 0$ is a certain tuning
parameter, and $\alpha \in (0,1)$ is introduced to control the compression errors. 
In algorithm \eqref{compressalgorithm}, 
only compressed information $Q\left(\cdot\right)$  is transmitted  from agent $j$ at iteration $k$, rather than the full-precision states. Specifically, agent $j$ computes the difference between its state $x_{j,k}$ and the auxiliary variable $h_{j,k}$, then applies the dynamic scaling factor $r_k$ and compression operator $Q(\cdot)$ to this difference, and transmits the compressed signal $Q(\cdot)$. Then, the receiver $i$ recovers $r_k Q(\cdot)$ to obtain  $\hat{x}_{j,k}$ as  in   \eqref{compress01}. 
The auxiliary variable $h_{j,k+1}$ is updated by \eqref{compress02}, where $\alpha$ controls the compression errors by smoothing the update of $h_{j,k+1}$. $\hat{x}_{W,i,k}$ is an auxiliary reconstruction variable that aggregates compressed information from neighbors to form the  neighborhood weighted consensus estimate in \eqref{compress11}. The variable update in \eqref{compress21} works as a backup copy for the neighboring information. Under the initial condition  $\tilde{{h}}_{i,0}=\sum_{j\in \mathcal{N}_i \cup \{i\}} w_{ij} {h}_{j,0}$, it holds that $\hat{{x}}_{{W},i,k}=\sum_{j\in \mathcal{N}_i \cup \{i\}} w_{ij}\hat{{x}}_{j,k}, \tilde{{h}}_{i,k}=\sum_{j\in \mathcal{N}_i \cup \{i\}} w_{ij} {h}_{j,k}$. This relation follows from mathematical induction, and the proof is omitted.

Algorithm \eqref{compressalgorithm} provides a communication-efficient implementation, where only compressed information $Q(\cdot)$ is transmitted. To facilitate convergence analysis, we derive an equivalent theoretical counterpart for a cleaner analysis, which is mathematically identical to Algorithm \eqref{compressalgorithm}. 
\begin{subequations} \label{compressalgorithmv}
	\begin{equation}  \label{equation1}
	x_{i,k+1}=x_{i,k}-\psi \sum_{j\in \mathcal{N}_{i}} w_{ij} (\hat{x}_{i,k}-\hat{x}_{j,k})+\tau (y_{i,k}-A_i z_{i,k}),
	\end{equation}
	\begin{equation}  \label{equation2}
	\tau y_{i,k+1}=\tau y_{i,k}-\psi \sum_{j\in \mathcal{N}_{i}} w_{ij}(\hat{x}_{i,k+1}-\hat{x}_{j,k+1}),
	\end{equation}
	\begin{equation}  \label{equation3}
	 z_{i,k+1}=z_{i,k}-\gamma \nabla f(z_{i,k}) +\gamma A_i^{T}(2 x_{i,k+1}-x_{i,k}),
	\end{equation}
\text{with}
		\begin{equation}  \label{compress1}
		\hat{x}_{i,k}=h_{i,k}+r_{k}Q\left(\frac{x_{i,k}-h_{i,k}}{r_{k}}\right),
	\end{equation}
	\begin{equation}  \label{compress2}
		h_{i,k+1}=(1-\alpha)h_{i,k}+\alpha \hat{x}_{i,k}.
	\end{equation}
\end{subequations}
By introducing $\mathbf{x}_{k}=\left[x_{1,k}^{T},...,x_{m,k}^{T}\right]^{T} \in \mathbb{R}^{nm}$, $\mathbf{y}_{k}=\left[y_{1,k}^{T},...,y_{m,k}^{T}\right]^{T}\in \mathbb{R}^{nm}$, $\mathbf{z}_{k}=\left[z_{1,k}^{T},...,z_{m,k}^{T}\right]^{T}\in \mathbb{R}^{nm}$, $\mathbf{h}_{k}=\left[h_{1,k}^{T},...,h_{m,k}^{T}\right]^{T}\in \mathbb{R}^{nm}$ and $\mathbf{A}=\text{blkdiag}(A_1,A_2,...,\\A_m) \in \mathbb{R}^{nm\times md}$, the algorithm in \eqref{compressalgorithm} is equivalent to the following compact form:
\begin{subequations} \label{compact}
\begin{equation}  \label{equation11}
	\mathbf{x}_{k+1}=\mathbf{x}_{k}-\psi (\mathbf{I}-\mathbf{W})\hat{\mathbf{x}}_{k} +\tau(\mathbf{y}_{k}-\mathbf{A}\mathbf{z}_{k}),
\end{equation}
\begin{equation}  \label{equation12}
	\tau \mathbf{y}_{k+1}=\tau \mathbf{y}_{k}-\psi (\mathbf{I}-\mathbf{W})\mathbf{\hat{x}}_{k+1},
\end{equation}
	\begin{equation}  \label{equation13}
	\mathbf{z}_{k+1}=\mathbf{z}_{k}-\gamma \nabla f(\mathbf{z}_{k}) +\gamma \mathbf{A}^{T}(2 \mathbf{x}_{k+1}-\mathbf{x}_{k}),
\end{equation}
\text{with}
	\begin{equation}  \label{equation14}
 \hat{\mathbf{x}}_{k}=\mathbf{h}_{k}+r_kQ(\frac{\mathbf{x}_{k}-\mathbf{h}_{k}}{r_{k}}),
\end{equation}
	\begin{equation}  \label{equation15}
	\mathbf{{h}}_{k+1}=(1-\alpha)\mathbf{h}_{k}+\alpha \hat{\mathbf{x}}_k,
\end{equation}
\end{subequations}
where $\mathbf{W} = W \otimes \mathbf{I}$. The iteration is initialized  by  $\mathbf{1}^{T}\mathbf{y}_0=\sum_{i=1}^{m}b_i, \mathbf{x}_0=\mathbf{z}_0=\hat{\mathbf{x}}_0=\mathbf{h}_0=\mathbf{0}$. If $\hat{\mathbf{x}}_{k}$ and $\hat{\mathbf{x}}_{k+1}$ are not compressed and $\psi=1$, algorithm \eqref{compact} will recover those in \eqref{DuSPA}.

\section{Convergence Analysis} \label{section4}
In this section, we provide convergence analysis of compressed algorithm \eqref{compact} under different compressors. 
We first establish several lemmas which are necessary for the subsequent analysis. Lemma \ref{equivalent} shows the equivalence between fixed-point of the steady-state form of \eqref{compact} and optimal solutions of  problem~\eqref{problem}. Lemma \ref{lemma2} guarantees the bijective transformation in the disagreement space. Lemma \ref{lemma3} provides the basis inequality for the convergence of the proposed algorithm.

\begin{lemma}  \label{equivalent}
	Suppose Assumptions \ref{existence}-\ref{matrix} hold. Under different compressors $Q_1(\cdot), Q_2(\cdot), Q_3(\cdot)$,  $\mathbf{z}^*$ is the optimal solution of problem \eqref{problem} if and only if the triple  $(\mathbf{x}^{\infty}, \mathbf{y}^{\infty}, \mathbf{z}^{\infty})$ is a fixed-point condition of the steady-state form of algorithm \eqref{compact}.
\end{lemma}
\begin{pf}
	First, if $\mathbf{z}^*$ is optimal for problem \eqref{problem}, according to the relation between problem \eqref{problem}, \eqref{dualproblem} and \eqref{equi}, the first-order optimality conditions guarantee the existence of
	 $\mathbf{x}^*$ such that  $\nabla f(\mathbf{z}^*)=\mathbf{A}^T\mathbf{x}^*,(\mathbf{I}-\mathbf{W})\mathbf{x}^*=\mathbf{0}$.  Meanwhile, there exists $\mathbf{y}^* = \mathbf{A} \mathbf{z}^*$, the tuple $(\mathbf{x}^*, \mathbf{y}^*, \mathbf{z}^*)$ satisfies fixed-point conditions of the steady-state form of algorithm \eqref{compact}.
	 
	 Conversely, if  $(\mathbf{x}^{\infty}, \mathbf{y}^{\infty}, \mathbf{z}^{\infty})$ is a fixed-point  of the steady-state form of \eqref{compact}, then the fixed-point conditions directly yield $\mathbf{A}^T\mathbf{x}^{\infty}=\nabla f(\mathbf{z}^{\infty})$ and  $\mathbf{y}^{\infty}=\mathbf{A}\mathbf{z}^{\infty}$. Furthermore, these conditions imply $\mathbf{h}^{\infty}=\mathbf{x}^{\infty}=\hat{\mathbf{x}}^{\infty}$, leading to $(\mathbf{I}-\mathbf{W})\hat{\mathbf{{x}}}^{\infty}=\mathbf{0}$ and $(\mathbf{I}-\mathbf{W})\mathbf{{{x}}}^{\infty}=\mathbf{0}$. These equalities are sufficient to establish that  $(\mathbf{x}^{\infty}, \mathbf{y}^{\infty}, \mathbf{z}^{\infty})$ satisfies the optimality conditions for problem  \eqref{problem}.  $\hfill\blacksquare$
\end{pf}

In this lemma, “a fixed point of the steady-state form” refers to the solution obtained when the algorithm’s
dynamic scaling factor $r_k\rightarrow 0$ as $k\rightarrow \infty$. This limiting
solution coincides with the optimal solution.

\begin{lemma} \label{lemma2}
\citep{xu2018bregman} Let $\mathbf{P}$ be a $m\times m$ matrix such that $null(\mathbf{P})=span\{\mathbf{1}\}$. Then, for each $\mathbf{y} \in span^{\perp}{\mathbf{1}}$, there exists a unique $\mathbf{y}'\in span^{\perp}\mathbf{1}$ such  that $\mathbf{y}=\mathbf{P}\mathbf{y}'$ and vice versa. 
\end{lemma}

\begin{lemma} \label{lemma3}
Let $\mathbf{s}_{k}=\left[\mathbf{x}_k^T, \mathbf{y}_k^T, \mathbf{z}_k^T\right]^{T}$ represent the sequence generated by the proposed compressed algorithm \eqref{compact}.  If $\gamma<\frac{\lambda_2}{\tau\rho_{A}}$ with $\lambda_2=\lambda_{\min}(\mathbf{W}), \rho_A=\lambda_{\max}(\mathbf{A}^T\mathbf{A})$,  the following inequality holds: 
\begin{equation}
	\begin{aligned}
\|\mathbf{s}_{k+1}-&\mathbf{s}^{*}\|_{\mathbf{M}}^2-\|\mathbf{s}_{k}-\mathbf{s}^{*}\|_{\mathbf{M}}^2+\|\Delta \mathbf{s}_{k+1}\|_{\mathbf{M}}^{2}\\
 & \le \tau\beta\left\|\Delta\mathbf{z}_{k+1}\right\|^{2}-\frac{2l_{f}L_{f}\tau}{L_{f}+l_{f}}\left\|\mathbf{z}_{k}-\mathbf{z}^{*}\right\|^{2} 
		\\
 &-\left(\frac{2\tau}{L_f+l_f}-\frac\tau\beta\right)\left\|\nabla f(\mathbf{z}_k)-\nabla f(\mathbf{z}^*)\right\|^2\\
 &+2\psi\langle(\boldsymbol{\varepsilon}_{k+1}-\boldsymbol{\varepsilon}_{k})(\mathbf{I}-\mathbf{W}),\mathbf{x}_{k+1}-\mathbf{x}^{*}\rangle\\
  &-2\tau\langle\mathbf{y}_{k+1}-\mathbf{y}^*,\boldsymbol{\varepsilon}_{k+1}\rangle,
		\end{aligned}
\end{equation}
where 
$\mathbf{M}=\begin{bmatrix}\mathbf{W}&\mathbf{0}&-\tau\mathbf{A}\\\mathbf{0}&\tau^2\mathbf{L}^{-1}/\psi &\mathbf{0}\\-\tau\mathbf{A}^{T}&\mathbf{0}&\frac{\tau}{\gamma}\mathbf{I}\end{bmatrix} \succ 0$, $\mathbf{L}= \mathbf{I}-(\mathbf{W}-\frac{\mathbf{1}\mathbf{1}^{T}}{m})$, $\beta$ is any positive number.

\end{lemma}
\begin{pf}
A detailed proof can be found in Appendix.  $\hfill\blacksquare$
\end{pf}

Based on Lemmas \ref{equivalent}-\ref{lemma3}, we aim to prove that the algorithm can achieve a linear convergence rate  under compressor $Q_1(\cdot)$ in the following theorem.

\begin{thm}  \label{thm1}
Under Assumptions \ref{existence}-\ref{matrix}, consider the compressor $Q_1(\cdot)$. Assume  $\gamma$ and $\tau$ satisfy $$\gamma <\min\left\{\frac{\lambda_{2}}{\tau}, \frac{\lambda_{2}}{4 \lambda_{2}\beta+\tau(4 \rho_{A}+1)},  \frac {2}{3} \frac {1}{L_{f}+l_{f}}\right\},$$  $\tau>\frac{L_f l_f}{\beta}$, where $\beta = L_f + l_f$.  Design the scaling factor follows $r_k^2=h\xi^k$ for some $h>0$. Then, algorithm \eqref{compact} converges linearly in expectation, i.e., there exist constants $\delta >0$  and $1-\delta < \xi < 1$ such that 
	\begin{equation}
 \mathbb{E}\left[\|\mathbf{s}_{k+1}-\mathbf{s}^*\|_{\mathbf{H}}^2\right]\le  (1-\delta)^{k+1}\mathbb{E}\left[\|\mathbf{s}_0-\mathbf{s}^*\|_{\mathbf{H}}^2\right]+c\xi^{k+1},
	\end{equation}
where $\mathbf{H}=\begin{bmatrix}\epsilon\mathbf{I}+\mathbf{W}&\mathbf{0}&-\tau\mathbf{A}\\\mathbf{0}&\tau^2(\epsilon\mathbf{I}+\mathbf{L}^{-1}/\psi )&\mathbf{0}\\-\tau\mathbf{A}^{T}&\mathbf{0}&\frac{\tau}{\gamma}\mathbf{I}\end{bmatrix} \succ 0$, $\epsilon=\gamma\tau$,  and the contraction parameter $\delta$ satisfies \[\delta < \min\left\{ \frac{\epsilon}{2\epsilon(1+\rho_A) + 2}, \\ \gamma\left( \frac{L_f l_f}{L_f + l_f} - \tau \right),\frac{\epsilon(1-\eta)}{2-\eta} \right\}.\]
\end{thm}
\begin{pf}    
A detailed proof can be found in Appendix.  $\hfill\blacksquare$
\end{pf}

Subsequently, we are going to prove that the algorithm \eqref{compact} achieves a linear convergence rate under compressor $Q_2(\cdot)$.

\begin{thm}   \label{thm2}
Under Assumptions  \ref{existence}-\ref{matrix}, consider the compressor $Q_2(\cdot)$. Let 
$$
\gamma < \min\left\{
\frac{2\lambda_2}{3\tau},\
\frac{2\lambda_2}{8\lambda_2\beta+\tau(8\rho_A+3)},\
\frac{2}{3}\left(
\frac{2}{L_f+l_f} - \frac{1}{\beta}
\right)
\right\},
$$
with $\tau>\frac{L_f l_f}{L_f + l_f}$. Define a Lyapunov function 
$\mathbb{E}\left[V_{k+1}\right]=\mathbb{E}\left[\|\mathbf{s}_{k+1}-\mathbf{s}^*\|_{\mathbf{H}}^2\right]+(a-d_1C)\mathbb{E}\left[\left\|\mathbf{x}_{k+1}-\mathbf{h}_{x}^{k+1}\right\| ^{2}\right]$, then there exists a constant $\nu \in (0,1)$ such that 
	\begin{equation}
		\mathbb{E}\left[V_{k+1}\right]\le (1-\nu) \mathbb{E}\left[V_{k}\right], ~~ \forall k\ge 0.
	\end{equation}
As a consequence, algorithm \eqref{compact} converges linearly with rate $(1-\nu)$, where the contraction coefficient $\nu$ satisfies $$1-\nu=\text{max}\left\{1-\delta, \frac{d_2C+ac_{x}}{a-d_1C}\right\},$$
where $d_1C<a<1$ with $C$ denoting the compression constant associated with  $Q_2(\cdot)$. 
\end{thm}

\begin{pf}    
A detailed proof can be found in Appendix.  $\hfill\blacksquare$
\end{pf}

Next, we will  extend our analysis to the more general case of biased compressors, demonstrating that algorithm \eqref{compact} achieves linear convergence under compressor $Q_3(\cdot)$. 

\begin{thm}   \label{thm3}
Under Assumptions \ref{existence}-\ref{matrix}, consider the compressor $Q_3(\cdot)$.  Let $\psi<\frac{1}{3\tau(1-\eta)}$ and $$ 2 < \gamma <\min\left\{\frac{\lambda_{2}}{\tau}, \frac{\lambda_{2}}{4 \lambda_{2}\beta+\tau(4 \rho_{A}+1)}, \frac{2}{3}\left(\frac{1}{L_{f}+l_{f}}\right) \right\},$$ with
$\tau>\frac{L_f l_f}{L_f + l_f}$.  Take the scaling factor $r_k^2=h\xi^k$ for some constant $h>0$, where  $1-\upsilon < \xi < 1$. Then, the sequence ${\mathbf{s}_k}$ generated by algorithm \eqref{compact} converges linearly. 
	\begin{equation}
		\|\mathbf{s}_{k+1}-\mathbf{s}^{*}\|_{\mathbf{H}}^{2}\le (1-\upsilon)^{k+1}\|\mathbf{s}_0-\mathbf{s}^*\|_{\mathbf{H}}^2+\varpi \xi^{k+1},
	\end{equation}
where $\upsilon$ satisfies $\upsilon < \min\left\{ \delta, \frac{(\epsilon - 2\tau)}{\epsilon+3\tau}\right\}.$
\end{thm}
\begin{pf}    
A detailed proof can be found in Appendix. $\hfill\blacksquare$
\end{pf}

Theorem \ref{thm1} and \ref{thm3} show  that the proposed algorithm \eqref{compact} achieves a linear convergence rate by setting $r_k^2=h \xi^{k}$ for constant $h >0$ and for some $\xi\in (0,1)$. Existing works (e.g., \cite{magnusson2020maintaining, xiong2022quantized}) achieve a linear convergence rate by designing specific time-varying uniform quantizers. Differently, our proposed algorithm \eqref{compact} attains linear convergence for a broad class of compressors.

\section{Numerical Simulation} \label{section5}

In this section, we use the numerical examples to  verify the effectiveness of the proposed algorithm under different compressors.

\subsection{Simulation Setting}

We simulate the effectiveness of the proposed algorithm using a quadratic optimization model. Each agent has a quadratic objective function, and the overall optimization objective is to minimize the total objective function $f(\mathbf{z})$ while satisfying the total load demand $\sum_{i=1}^{m}b_i$. The problem can be formulated as $f(\mathbf{z})=\sum_{i=1}^{m}\alpha_iz_i^2+\beta_iz_i+\gamma_i$. 
where $\alpha_i,\beta_i,\gamma_i>0$ are the coefficients and $z_i$ is the decision variable of agent $i$. The parameters are adopted from \citep{xu2018dual} and restated in Table \ref{parameters}. The coupling matrix elements we set $A_i=1$ and the load demand $\sum_{i}^{m}b_i=0$.
 
\begin{table}[!h]
	\renewcommand{\arraystretch}{1.0}
	\caption{System~~Parameters} 
	\label{parameters} 
	\centering 
	\begin{tabular}{|c|c|c|c|} 
		\hline 
		Bus & $\alpha_i$ (\$/MW$^2$) & $\beta_i$ (\$/MW) & $\gamma_i$ (\$)\\ 
		\hline 
		1   & 0.04            & 2.0        & 0       \\ 
		\hline 
		2   & 0.03            & 3.0        & 0         \\ 
		\hline 
		3   & 0.035           & 4.0        & 0        \\ 
		\hline 
		6   & 0.03            & 4.0        & 0     \\ 
		\hline 
		8   & 0.04            & 2.5        & 0       \\ 
		\hline 
	\end{tabular}
\end{table}

In the following cases, we use $x$ to denote the input of compressors $Q(\cdot)$ for brevity.

\textbf {Case 1:}
Under compressor $Q_1(\cdot)$, we use the unbiased probabilistic quantizer \citep{yuan2012distributed} to compress the transmitted signal. 
\begin{equation}
{Q}_1 (x) = 
\begin{cases} 
	\lfloor x \rfloor_p & \text{with probability } (\lceil x \rceil_p - x) \Delta_p \\
	\lceil x \rceil_p & \text{with probability } (x - \lfloor x \rfloor_p) \Delta_p, \nonumber 
\end{cases}
\end{equation}
where $\lfloor x \rfloor_p$ and $\lceil x \rceil_p$ denote the operation of rounding down and up $x$ to the nearest integer multiples of $1/\Delta_p$, respectively, and $\Delta_p$ is some positive integer.

\textbf {Case 2:}
Under compressor $Q_2(\cdot)$, we use the unbiased  $b_{Q_2}$-bits quantization method with $\infty$-norm in \citep{liao2022compressed}
$
{Q}_2(x) = \left( \|x\|_{\infty} \cdot 2^{-(b_{Q_2}-1)} \cdot \text{sign}(x) \right) \odot \left\lfloor \frac{2^{(b_{Q_2}-1)} \cdot |x|}{\|x\|_{\infty}} + \mu \right\rfloor$ to compress the transmitted signal, where $\odot$ is the Hadamard product, and $\mu$ is a random vector uniformly distributed in $[0,1]^{d}$. 

\textbf {Case 3:}
Under compressor $Q_3(\cdot)$, we use the biased deterministic  truncation  quantizer ${Q}_{3}(x)=	\lfloor x \rfloor$ in \citep{el2016design} to compress the transmitted signal,  which rounds the value to the
nearest lower integer.

\subsection{Simulation Results}

This part provides the convergence results along with a corresponding analysis for three specific quantizers.

\begin{figure}[h]
	\centering
	\begin{minipage}{4cm}
		\centering
		\includegraphics[width=4.5cm]{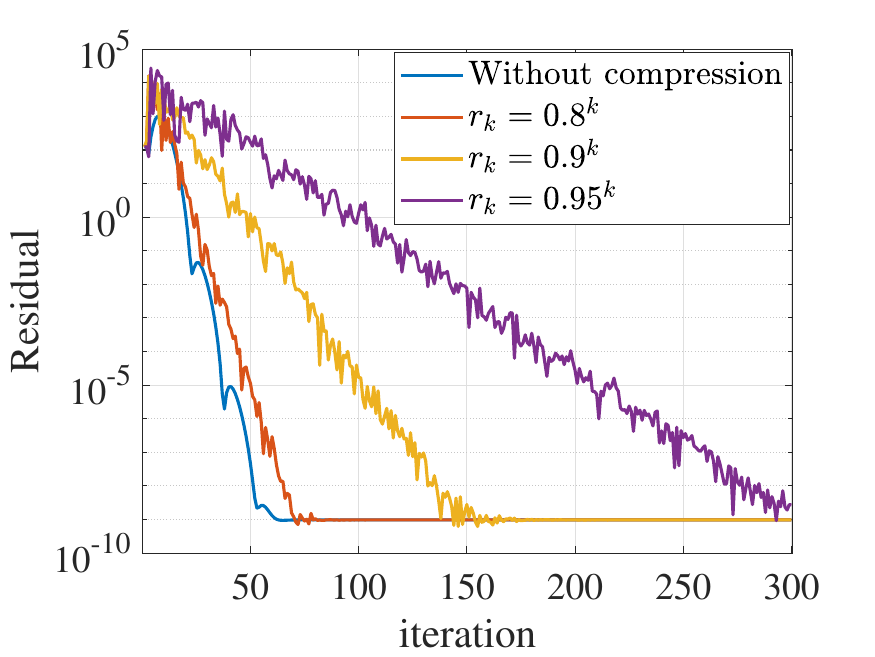}
		\caption{Residuals under different $r_k$.}
	  \label{compression11}
	\end{minipage}
	\quad
	\begin{minipage}{4cm}
		\centering
		\includegraphics[width=4.5cm]{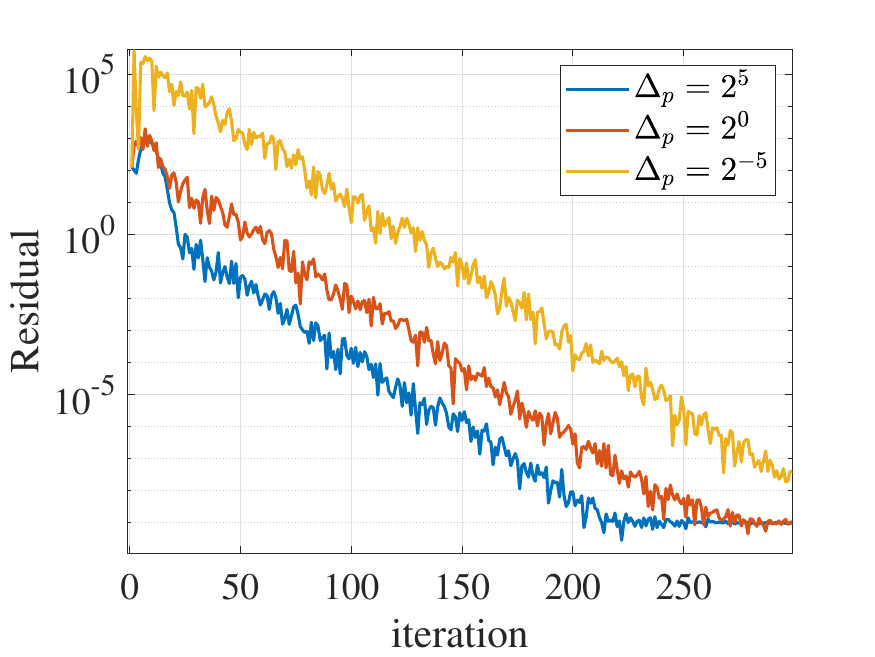}
		\caption{Residuals under different $\Delta_p$.}
    	\label{compression12}
	\end{minipage}
\end{figure}

\textbf {Case 1:}
Under quantizer $Q_{1}(x)$, the simulation results are given in Fig.\ref{compression11} and Fig.\ref{compression12}.   Fig.\ref{compression11} demonstrates how the dynamic scaling factor regulates the performance of the compressed algorithm. It indicates that as $r_k$ increases, the convergence performance of the algorithm deteriorates. This occurs because increasing $r_k$ weakens the attenuation effect on the absolute compression error. 

Fig.\ref{compression12} indicates that under the same $r_k$ condition, a larger $\Delta_p$ leads to a  worse algorithm convergence performance. This is because a larger $\Delta_p$ leads to a greater interval between discrete values, requiring fewer bits to distinguish them. However, this will result in a larger quantization error, thus more severely impacting the algorithm's convergence performance.

\begin{figure}[h]
	\centering
	\begin{minipage}{4cm}
		\centering
		\includegraphics[width=4.5cm]{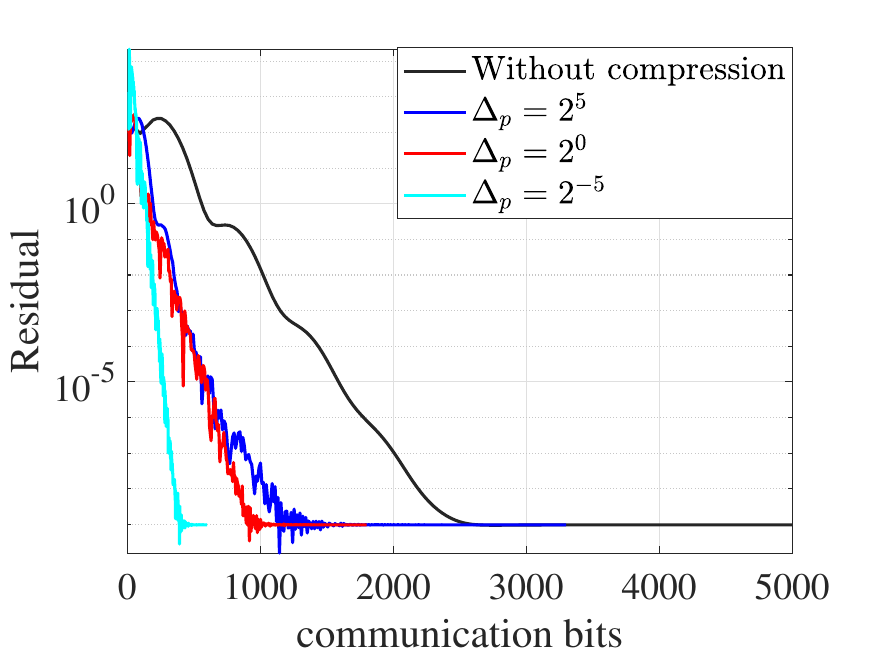}
		\caption{Residuals vs. Communication bits}
		\label{compression13}
	\end{minipage}
	\quad
	\begin{minipage}{4cm}
		\centering
		\includegraphics[width=4.5cm]{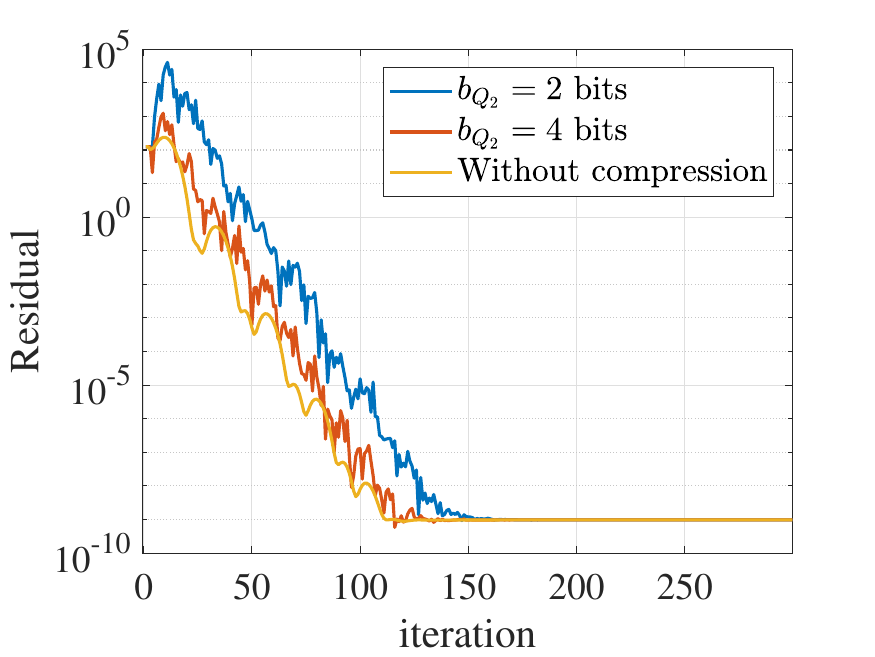}
		\caption{Residuals under different $b_{Q_2}$.}
		\label{compression2}
	\end{minipage}
\end{figure}

Under quantizer ${Q}_1 (x)$, the transmitted bits of input $x$ can be denoted as $b_{Q_1}=\lceil \log_2(\lfloor (2C_{max})\Delta_p\rfloor +1)\rceil$. $C_{max}$  is the uniform bound of the scaled compressor input $|\frac{x_k-h_k}{r_k}|$,  which is $\mathcal{O}(1)$ due to matched decay rates between the difference $(x_k-h_k)$ and the dynamic scaling factor $r_k$, and is obtained by offline pre-running the algorithm. In the absence of quantization, our approach adopts the prevailing 32-bit data format, as utilized in traditional digital communication systems \citep{roberts1987digital}. From Fig.\ref{compression13}, a key finding is that the communication cost, measured in total bits transmitted to achieve a specific residual error, is monotonically increasing with the quantization parameter $\Delta_p$.  However, the compressed version requires fewer transmitted bits for
transmission compared to the uncompressed scheme, while maintaining the same level of accuracy.

\textbf {Case 2:}
Under quantizer $Q_{2}(x)$, Fig.\ref{compression2} presents a comparison of the convergence performance under different quantization bits versus the uncompressed algorithm. We can observe that compressed communication impairs the convergence rate of the algorithm. However, as the number of transmitted bits increases, the communication conditions improve, and the convergence rate of the algorithm also improves.

\begin{figure}[h]
	\centering
	\begin{minipage}{4cm}
		\centering
		\includegraphics[width=4.5cm]{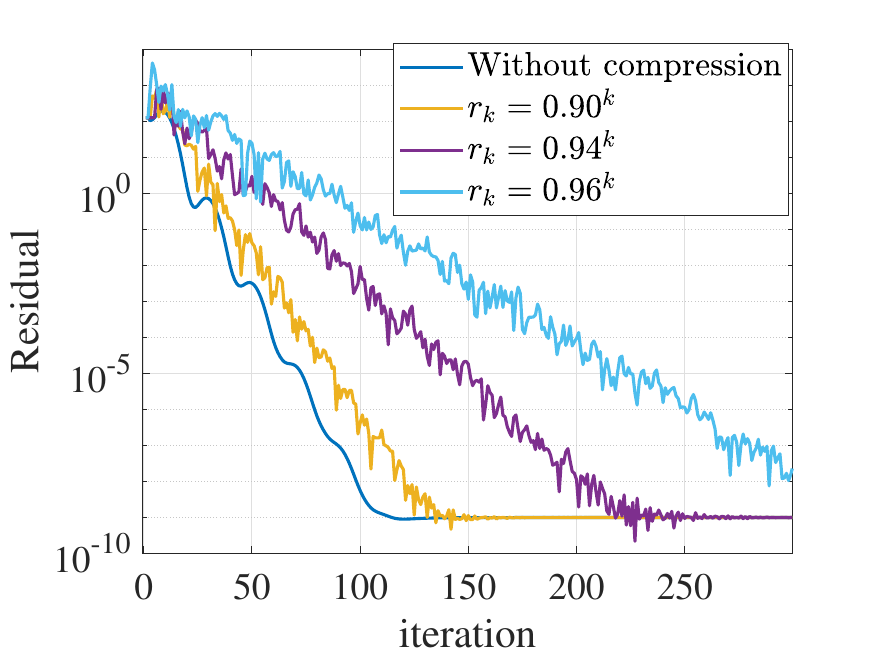}
		\caption{Residuals under different $r_k$.}
		\label{compression3}
	\end{minipage}
	\quad
	\begin{minipage}{4cm}
		\centering
		\includegraphics[width=4.5cm]{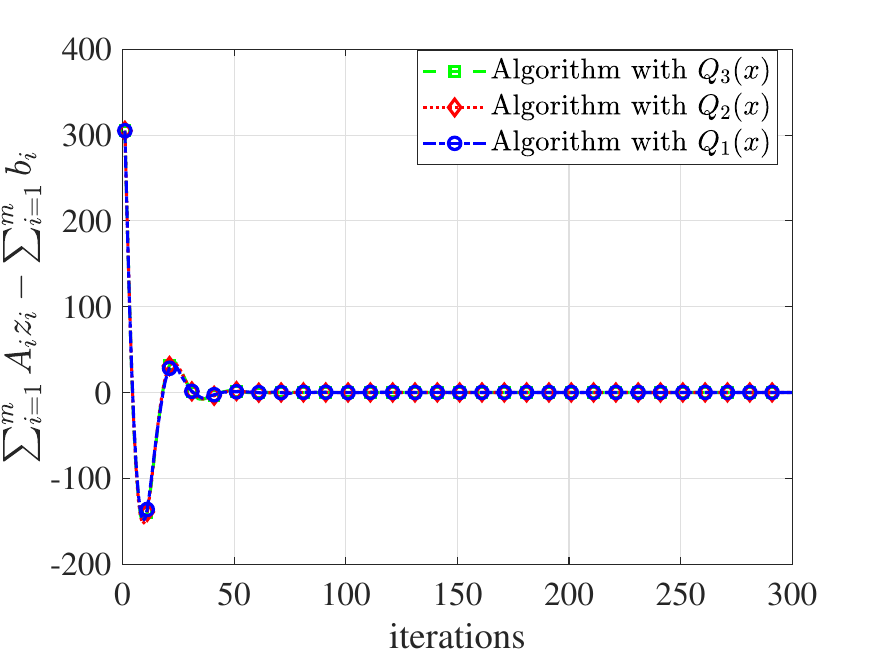}
		\caption{Equality constraint violations.}
		\label{Residual}
	\end{minipage}
\end{figure}

\textbf {Case 3:}
Under quantizer $Q_{3}(x)$, Fig.\ref{compression3} shows that the smaller $r_k$ is, the better the suppression effect on the biased quantization error, thereby making the convergence performance closer to that of lossless transmission.

Finally, we set $r_k = 0.96^k, \gamma = 3, \tau=0.03$, and the number of transmission bits $b_{Q_2}=2$. Fig.\ref{Residual} shows that the algorithm's equality constraints are satisfied under different quantizers. This validates the effectiveness of the proposed algorithm with compression in handling coupled equality constraints.

\section{Conclusion} \label{section6}

In this paper, we propose a distributed optimization algorithm with compressed communication to address problems involving coupled equality constraints.  By introducing a difference compression technique and dynamic scaling factors, the algorithm achieves linear convergence under different types of compressors while satisfying coupled equality constraints. The results confirm that the same convergence accuracy can be achieved with reduced communication overhead. An important direction for future work is to generalize our analysis, to validate the effectiveness of the algorithm under directed graphs and more complex constraint scenarios.

\bibliography{ifacconf}    

@inproceedings{Qian2025,
	title={Quantized Distributed Dual Algorithm for Resource Allocation with Gradient Tracking over Unbalanced Networks},
	author={Qian, Kun and Li, Yun-Long and Liu, Xiao-Kang},
	booktitle={44th Chinese Control Conference (CCC)},
	year={2025},
	pages={2100--2105},
}

@Inproceedings{Duan2024ASCC,
	author={Duan, Yuzhu and Yang, Ziwen and Zhu, Shanying},
	booktitle={14th Asian Control Conference (ASCC)}, 
	title={A Communication-Efficient Distributed Optimization Algorithm with Linear Convergence over Directed Graphs}, 
	year={2024},
	pages={1080--1085},
}

@article{guo2020adaptive,
	title={An adaptive wireless virtual reality framework in future wireless networks: A distributed learning approach},
	author={Guo, Fengxian and Yu, F Richard and Zhang, Heli and Ji, Hong and Leung, Victor CM and Li, Xi},
	journal={IEEE Transactions on Vehicular Technology},
	volume={69},
	number={8},
	pages={8514--8528},
	year={2020},
	publisher={IEEE}
}

@article{lee2021fast,
	title={Fast and scalable distributed consensus over wireless large-scale internet of things network},
	author={Lee, Hojung and Seo, Hyowoon and Choi, Wan},
	journal={IEEE Internet of Things Journal},
	volume={9},
	number={11},
	pages={7916--7930},
	year={2021},
	publisher={IEEE}
}

@article{ding2021differentially,
	title={Differentially private distributed resource allocation via deviation tracking},
	author={Ding, Tie and Zhu, Shanying and Chen, Cailian and Xu, Jinming and Guan, Xinping},
	journal={IEEE Transactions on Signal and Information Processing over Networks},
	volume={7},
	pages={222--235},
	year={2021},
	publisher={IEEE}
}

@article{erol2012multi,
	title={A multi-agent based approach to dynamic scheduling of machines and automated guided vehicles in manufacturing systems},
	author={Erol, Rizvan and Sahin, Cenk and Baykasoglu, Adil and Kaplanoglu, Vahit},
	journal={Applied Soft Computing},
	volume={12},
	number={6},
	pages={1720--1732},
	year={2012},
	publisher={Elsevier}
}

@article{nishi2005distributed,
	title={Distributed route planning for multiple mobile robots using an augmented Lagrangian decomposition and coordination technique},
	author={Nishi, Tatsushi and Ando, Masakazu and Konishi, Masami},
	journal={IEEE Transactions on Robotics},
	volume={21},
	number={6},
	pages={1191--1200},
	year={2005},
	publisher={IEEE}
}

@article{zhang2018framework,
	title={A framework for smart production-logistics systems based on CPS and industrial IoT},
	author={Zhang, Yingfeng and Guo, Zhengang and Lv, Jingxiang and Liu, Ying},
	journal={IEEE Transactions on Industrial Informatics},
	volume={14},
	number={9},
	pages={4019--4032},
	year={2018},
	publisher={IEEE}
}

@inproceedings{karimireddy2019error,
	title={Error feedback fixes signsgd and other gradient compression schemes},
	author={Karimireddy, Sai Praneeth and Rebjock, Quentin and Stich, Sebastian and Jaggi, Martin},
	booktitle={36th International Conference on Machine Learning (ICML)},
	pages={3252--3261},
	year={2019},
}

@inproceedings{alistarh2017qsgd,
	title={{QSGD}: Communication-efficient SGD via gradient quantization and encoding},
	author={Alistarh, Dan and Grubic, Demjan and Li, Jerry and Tomioka, Ryota and Vojnovic, Milan},
	booktitle={31th Advances in Neural Information Processing Systems (NeurIPS)},
	volume={30},
	year={2017},
}

@article{zhu2018mitigating,
	title={Mitigating quantization effects on distributed sensor fusion: A least squares approach},
	author={Zhu, Shanying and Chen, Cailian and Xu, Jinming and Guan, Xinping and Xie, Lihua and Johansson, Karl Henrik},
	journal={IEEE Transactions on Signal Processing},
	volume={66},
	number={13},
	pages={3459--3474},
	year={2018},
	publisher={IEEE}
}

@inproceedings{nedic2008distributed,
	title={Distributed subgradient methods and quantization effects},
	author={Nedic, Angelia and Olshevsky, Alex and Ozdaglar, Asuman and Tsitsiklis, John N},
	booktitle={47th IEEE Conference on Decision and Control (CDC)},
	pages={4177--4184},
	year={2008},
}

@article{aysal2008distributed,
	title={Distributed average consensus with dithered quantization},
	author={Aysal, Tuncer Can and Coates, Mark J and Rabbat, Michael G},
	journal={IEEE Transactions on Signal Processing},
	volume={56},
	number={10},
	pages={4905--4918},
	year={2008},
	publisher={IEEE}
}

@inproceedings{tang2019doublesqueeze,
	title={Doublesqueeze: Parallel stochastic gradient descent with double-pass error-compensated compression},
	author={Tang, Hanlin and Yu, Chen and Lian, Xiangru and Zhang, Tong and Liu, Ji},
	booktitle={36th International Conference on Machine Learning (ICML)},
	pages={6155--6165},
	year={2019},
}

@inproceedings{tang2018communication,
	title={Communication compression for decentralized training},
	author={Tang, Hanlin and Gan, Shaoduo and Zhang, Ce and Zhang, Tong and Liu, Ji},
	booktitle={32nd Advances in Neural Information Processing Systems (NeurIPS)},
	volume={31},
	pages={7652--7662},
	year={2018},
}

@inproceedings{koloskova2019decentralized,
	title={Decentralized stochastic optimization and gossip algorithms with compressed communication},
	author={Koloskova, Anastasia and Stich, Sebastian and Jaggi, Martin},
	booktitle={36th International Conference on Machine Learning (ICML)},
	pages={3478--3487},
	year={2019},
}

@inproceedings{koloskova2019decentralized1,
	title={Decentralized deep learning with arbitrary communication compression},
	author={Koloskova, Anastasia and Lin, Tao and Stich, Sebastian U and Jaggi, Martin},
	booktitle={7th International Conference on Learning Representations (ICLR)},
	year={2019},
}

@article{song2022compressed,
	title={Compressed gradient tracking for decentralized optimization over general directed networks},
	author={Song, Zhuoqing and Shi, Lei and Pu, Shi and Yan, Ming},
	journal={IEEE Transactions on Signal Processing},
	volume={70},
	pages={1775--1787},
	year={2022},
	publisher={IEEE}
}

@article{xiong2022quantized,
	title={Quantized distributed gradient tracking algorithm with linear convergence in directed networks},
	author={Xiong, Yongyang and Wu, Ligang and You, Keyou and Xie, Lihua},
	journal={IEEE Transactions on Automatic Control},
	volume={68},
	number={9},
	pages={5638--5645},
	year={2023},
	publisher={IEEE}
}

@article{magnusson2020maintaining,
	title={On maintaining linear convergence of distributed learning and optimization under limited communication},
	author={Magn{\'u}sson, Sindri and Shokri-Ghadikolaei, Hossein and Li, Na},
	journal={IEEE Transactions on Signal Processing},
	volume={68},
	pages={6101--6116},
	year={2020},
	publisher={IEEE}
}

@article{shi2024distributed,
	title={Distributed economic dispatch algorithm with quantized communication mechanism},
	author={Shi, Xiasheng and Sun, Changyin and Mu, Chaoxu},
	journal={IEEE Transactions on Automation Science and Engineering},
	volume={22},
	pages={8618--8629},
	year={2024},
	publisher={IEEE}
}

@article{liao2022compressed,
	title={A compressed gradient tracking method for decentralized optimization with linear convergence},
	author={Liao, Yiwei and Li, Zhuorui and Huang, Kun and Pu, Shi},
	journal={IEEE Transactions on Automatic Control},
	volume={67},
	number={10},
	pages={5622--5629},
	year={2022},
	publisher={IEEE}
}

@inproceedings{liulinear,
	title={Linear Convergent Decentralized Optimization with Compression},
	author={Liu, Xiaorui and Li, Yao and Wang, Rongrong and Tang, Jiliang and Yan, Ming},
	booktitle={9th International Conference on Learning Representations (ICLR)},
    year={2021},
}

@article{wang2020dual,
	title={A dual splitting method for distributed economic dispatch in multi-energy systems},
	author={Wang, Zhibin and Xu, Jinming and Zhu, Shanying and Chen, Cailian},
	journal={IFAC-PapersOnLine},
	volume={53},
	number={2},
	pages={12566--12571},
	year={2020},
	publisher={Elsevier}
}

@article{ren2026distributed,
	title={Distributed Optimization by Network Flows With Spatio-Temporal Compression},
	author={Ren, Zihao and Wang, Lei and Yi, Xinlei and Wang, Xi and Yuan, Deming and Yang, Tao and Wu, Zhengguang and Shi, Guodong},
	journal={IEEE Transactions on Automatic Control},
	year={2026},
	pages={1--16},
	publisher={IEEE}
}

@article{mishchenko2025distributed,
	title={Distributed learning with compressed gradient differences},
	author={Mishchenko, Konstantin and Gorbunov, Eduard and Tak{\'a}{\v{c}}, Martin and Richt{\'a}rik, Peter},
	journal={Optimization Methods and Software},
	volume={40},
	number={5},
	pages={1181--1196},
	year={2025},
}

@article{xu2018bregman,
	title={A Bregman splitting scheme for distributed optimization over networks},
	author={Xu, Jinming and Zhu, Shanying and Soh, Yeng Chai and Xie, Lihua},
	journal={IEEE Transactions on Automatic Control},
	volume={63},
	number={11},
	pages={3809--3824},
	year={2018},
	publisher={IEEE}
}

@article{xu2018dual,
	title={A dual splitting approach for distributed resource allocation with regularization},
	author={Xu, Jinming and Zhu, Shanying and Soh, Yeng Chai and Xie, Lihua},
	journal={IEEE Transactions on Control of Network Systems},
	volume={6},
	number={1},
	pages={403--414},
	year={2019},
	publisher={IEEE}
}

@article{yuan2012distributed,
	title={Distributed dual averaging method for multi-agent optimization with quantized communication},
	author={Yuan, Deming and Xu, Shengyuan and Zhao, Huanyu and Rong, Lina},
	journal={Systems \& Control Letters},
	volume={61},
	number={11},
	pages={1053--1061},
	year={2012},
	publisher={Elsevier}
}

@article{el2016design,
	title={Design and analysis of distributed averaging with quantized communication},
	author={El Chamie, Mahmoud and Liu, Ji and Ba{\c{s}}ar, Tamer},
	journal={IEEE Transactions on Automatic Control},
	volume={61},
	number={12},
	pages={3870--3884},
	year={2016},
	publisher={IEEE}
}

@book{roberts1987digital,
	title={Digital Signal Processing},
	author={Roberts, Richard A and Mullis, Clifford T},
	year={1987},
	publisher={Addison-Wesley Longman Publishing Co., Inc.}
}

\newpage 
\appendix
\section{Proof of Lemma 3}    
According to (\ref{equation11}) and (\ref{equation12}), we can obtain,
\begin{equation}  \label{deltax}
	\begin{aligned}
		\mathbf{W}\Delta \mathbf{x}_{k+1}=\tau\mathbf{y}_{k+1}-\tau\mathbf{A}\mathbf{z}_{k}+\psi(\boldsymbol{\varepsilon}_{k+1}-\boldsymbol{\varepsilon}_{k})(\mathbf{I}-\mathbf{W}), 
	\end{aligned}
\end{equation}
it can be further obtained that,
\begin{equation} \label{xx}
	\begin{split}
		&\mathbf{W}\langle\mathbf{x}_{k+1}-\mathbf{x}_{k},\mathbf{x}_{k+1}-\mathbf{x}^{*}\rangle\\
	&=	\langle\tau\mathbf{y}_{k+1}-\tau\mathbf{y}^{*}+\tau\mathbf{A}(\mathbf{z}^{*}-\tau\mathbf{z}_{k})\\
		&\quad +\psi(\boldsymbol{\varepsilon}_{k+1}-\boldsymbol{\varepsilon}_{k})(\mathbf{I}-\mathbf{W}), \mathbf{x}_{k+1}-\mathbf{x}^{*}\rangle\\
		&=\tau\langle\mathbf{y}_{k+1}-\mathbf{y}^{*},\mathbf{x}_{k+1}-\mathbf{x}^{*}\rangle-\tau\langle\mathbf{A}(\mathbf{z}_{k}-\mathbf{z}^{*}),\\
		&\quad \mathbf{x}_{k+1}-\mathbf{x}^{*}\rangle+\psi\langle(\boldsymbol{\varepsilon}_{k+1}-\boldsymbol{\varepsilon}_{k})(\mathbf{I}-\mathbf{W}),\mathbf{x}_{k+1}-\mathbf{x}^{*}\rangle\\
		&=\tau\langle\mathbf{y}_{k+1}-\mathbf{y}^{*},\mathbf{x}_{k+1}-\mathbf{x}^{*}\rangle-\tau\mathbf{A}\langle\mathbf{z}_{k}+\mathbf{z}_{k+1}\\
		&\quad-\mathbf{z}_{k+1}-\mathbf{z}^{*},\mathbf{x}_{k+1}-\mathbf{x}^{*}\rangle+\psi\langle(\boldsymbol{\varepsilon}_{k+1}\\
		&\quad-\boldsymbol{\varepsilon}_{k})(\mathbf{I}-\mathbf{W}),\mathbf{x}_{k+1}-\mathbf{x}^{*}\rangle\\
		&=\tau\langle\mathbf{y}_{k+1}-\mathbf{y}^{*},\mathbf{x}_{k+1}-\mathbf{x}^{*}\rangle\\
		&\quad+\tau\langle\mathbf{A}(\mathbf{z}_{k+1}-\mathbf{z}_{k}),\mathbf{x}_{k+1}-\mathbf{x}^{*}\rangle\\
		&\quad-\tau\langle\mathbf{A}(\mathbf{z}_{k+1}-\mathbf{z}^{*}),\mathbf{x}_{k+1}-\mathbf{x}^{*}\rangle\\
		&\quad+\psi\langle(\boldsymbol{\varepsilon}_{k+1}-\boldsymbol{\varepsilon}_{k})(\mathbf{I}-\mathbf{W}),\mathbf{x}_{k+1}-\mathbf{x}^{*}\rangle,
	\end{split}
\end{equation}
if we let $\mathbf{y}_{k+1}-\mathbf{y}^{*}=(\mathbf{I}-\mathbf{W})(\mathbf{y}_{k+1}^{\prime}-\mathbf{y}^{*\prime})$ and define $\mathbf{L}= \mathbf{I}-(\mathbf{W}-\frac{\mathbf{1}\mathbf{1}^{T}}{m})$, then we can obtain
\begin{equation}
	\begin{aligned}&\langle\mathbf{y}_{k+1}-\mathbf{y}^{*},\mathbf{x}_{k+1}-\mathbf{x}^{*}\rangle\\
		&=\langle\mathbf{y}_{k+1}^{\prime}-\mathbf{y}^{*\prime},(\mathbf{I}-\mathbf{W})(\mathbf{x}_{k+1}-\mathbf{x}^{*})\rangle\\
		&=\langle\mathbf{y}_{k+1}^{\prime}-\mathbf{y}^{*\prime},(\mathbf{I}-\mathbf{W})(\mathbf{x}_{k+1}-\hat{\mathbf{x}}_{k+1}\\
		&\quad+\hat{\mathbf{x}}_{k+1}-\mathbf{x}^{*})\rangle\\
		&=\langle\mathbf{y}_{k+1}^{\prime}-\mathbf{y}^{*\prime},(\mathbf{I}-\mathbf{W})(\mathbf{x}_{k+1}-\hat{\mathbf{x}}_{k+1})\rangle\\
		&+\langle\mathbf{y}_{k+1}^{\prime}-\mathbf{y}^{*\prime},(\mathbf{I}-\mathbf{W})(\hat{\mathbf{x}}_{k+1}-\mathbf{x}^{*})\rangle\\
		&={-\tau/\psi\langle\mathbf{y}_{k+1}-\mathbf{y}^{*},\mathbf{y}_{k+1}-\mathbf{y}_{k}\rangle_{\mathbf{L}^{-1}}}\\
		&\quad-\langle\mathbf{y}_{k+1}-\mathbf{y}^{*},\boldsymbol{\varepsilon}_{k+1}\rangle.
	\end{aligned}
\end{equation}

According to (\ref{equation11}) and (\ref{equation13}), we can obtain that
\begin{equation}  \label{xz}
	\begin{aligned}
		&\langle\mathbf{A}(\mathbf{z}_{k+1}-\mathbf{z}^{*}),\mathbf{x}_{k+1}-\mathbf{x}^{*}\rangle\\
		&=\langle\mathbf{y}_{k+1}-\mathbf{y}^{*},\mathbf{x}_{k+1}-\mathbf{x}^{*}\rangle\\
		&\quad+\langle\mathbf{A}(\mathbf{z}_{k+1}-\mathbf{z}_{k}),\mathbf{x}_{k+1}-\mathbf{x}^{*}\rangle\\
		&\quad+\frac{1}{\tau}\langle\psi(\boldsymbol{\varepsilon}_{k+1}-\boldsymbol{\varepsilon}_{k})(\mathbf{I}-\mathbf{W}),\mathbf{x}_{k+1}-\mathbf{x}^{*}\rangle\\
		&\quad-\frac{1}{\tau}\mathbf{W}\langle\mathbf{x}_{k+1}-\mathbf{x}_{k},\mathbf{x}_{k+1}-\mathbf{x}^{*}\rangle.
	\end{aligned}
\end{equation}

Substitute (\ref{xx}) and (\ref{xz}) into (\ref{deltax}), it can be deduced that,
\begin{equation} \label{zz}
	\begin{aligned}
		&{2\langle\mathbf{z}_{k+1}-\mathbf{z}^{*},\mathbf{z}_{k+1}-\mathbf{z}_{k}\rangle}-2\gamma\langle\mathbf{A}(\mathbf{z}_{k+1}-\mathbf{z}^{*}),\\
		&\quad{\mathbf{x}}_{k+1}-\mathbf{x}^{*}\rangle-2\gamma\langle\mathbf{A}(\mathbf{z}_{k+1}-\mathbf{z}^{*}),{\mathbf{x}}_{k+1}-\mathbf{x}_{k}\rangle\\
		&=2\langle\mathbf{z}_{k+1}-\mathbf{z}^*,\mathbf{z}_{k+1}-\mathbf{z}_k\rangle-2\gamma\langle\mathbf{A}(\mathbf{z}_{k+1}-\mathbf{z}^*),\\
		&\quad\mathbf{x}_{k+1}-\mathbf{x}_k\rangle-\frac{2\gamma \mathbf{W}}\tau\langle\mathbf{x}_{k+1}-\mathbf{x}^*,\mathbf{x}_k-\mathbf{x}_{k+1}\rangle \\
		&\quad-2\gamma\langle\mathbf{y}_{k+1}-\mathbf{y}^*,\mathbf{x}_{k+1}-\mathbf{x}^*\rangle-2\gamma\langle\mathbf{A}(\mathbf{z}_{k+1}-\mathbf{z}_k),\\
		&\quad\mathbf{x}_{k+1} -\mathbf{x}^*\rangle -\frac{2\gamma}\tau\langle(\psi \boldsymbol{\varepsilon}_{k+1}-\boldsymbol{\varepsilon}_{k})(\mathbf{I}-\mathbf{W}),\mathbf{x}_{k+1}-\mathbf{x}^{*}\rangle\\
		&=2\langle\mathbf{z}_{k+1}-\mathbf{z}^*,\mathbf{z}_{k+1}-\mathbf{z}_k\rangle+\frac{2\gamma \mathbf{W}}\tau\langle\mathbf{x}_{k+1} \\
		&\quad-\mathbf{x}^*,\mathbf{x}_{k+1}-\mathbf{x}_k\rangle -2\gamma\langle\mathbf{A}(\mathbf{z}_{k+1}-\mathbf{z}^*),\mathbf{x}_{k+1}-\mathbf{x}_k\rangle\\
		&\quad-2\gamma\langle\mathbf{A}(\mathbf{z}_{k+1}-\mathbf{z}_k),\mathbf{x}_{k+1}-\mathbf{x}^*\rangle \\
		&\quad+2\gamma\tau/\psi \langle\mathbf{y}_{k+1}-\mathbf{y}^*,\mathbf{y}_{k+1}-\mathbf{y}_k\rangle_{\mathbf{L}^{-1}}\\
		&\quad+2\gamma\langle\mathbf{y}_{k+1}-\mathbf{y}^*,\boldsymbol{\varepsilon}_{k+1}\rangle\nonumber
			\end{aligned}
	\end{equation}
		\begin{equation} 
			\begin{aligned}
		&-\frac{2\gamma}\tau\psi \langle(\boldsymbol{\varepsilon}_{k+1}-\boldsymbol{\varepsilon}_{k})(\mathbf{I}-\mathbf{W}),\mathbf{x}_{k+1}-\mathbf{x}^{*}\rangle.
	\end{aligned}
\end{equation}
Recalling that $\Delta \mathbf{z}_{k+1}+\gamma \mathbf{A}^{T}(2\mathbf{x}_{k+1}-\mathbf{x}_k)-\gamma \nabla f(\mathbf{z}_k)= 0 $ from \eqref{equation13}, together with $\gamma \mathbf{A}\mathbf{x}^{*}-\gamma \nabla f(\mathbf{z}^{*})=0 $ from the optimality condition, and knowing that $\partial g$ is maximally monotone we have $\langle \Delta \mathbf{z}_{k+1}+\gamma (2\mathbf{x}_{k+1}-\mathbf{x}^*)-\gamma(\nabla f(\mathbf{z}_k)-\nabla f(\mathbf{z}^*)), \mathbf{z}_{k+1}-\mathbf{z}^*=0$. Substituting the above relation into \eqref{zz} and then multiplying both sides by  $\frac{\tau}{\gamma}$, and we can obtain:

\begin{equation} \label{tgamma}
	\begin{aligned}
		&2\frac{\tau}{\gamma}\langle\mathbf{z}_{k+1}-\mathbf{z}^*,\mathbf{z}_{k+1}-\mathbf{z}_k\rangle+2\mathbf{W}\langle\mathbf{x}_{k+1}-\mathbf{x}^*,\\
		&\quad\mathbf{x}_{k+1}-\mathbf{x}_k\rangle+2\tau^2/\psi\langle\mathbf{y}_{k+1}-\mathbf{y}^*,\\
		&\quad\mathbf{y}_{k+1}-\mathbf{y}_k\rangle_{(\mathbf{I}-\mathbf{W})^{-1}} \\
		&\quad-2\tau\langle\mathbf{A}(\mathbf{z}_{k+1}-\mathbf{z}^*),\mathbf{x}_{k+1}-\mathbf{x}_k\rangle\\
		&\quad-2\tau\langle\mathbf{A}(\mathbf{z}_{k+1}-\mathbf{z}_k),\mathbf{x}_{k+1}-\mathbf{x}^*\rangle \\
		&= -2\tau\langle\mathbf{z}_{k+1}-\mathbf{z}^{*}, \nabla f(\mathbf{z}_{k})-\nabla  f(\mathbf{z}^{*})\rangle\\
		&\quad-2\tau\langle\mathbf{y}_{k+1}-\mathbf{y}^*,\boldsymbol{\varepsilon}_{k+1}\rangle\\
		&\quad+2\psi \langle(\boldsymbol{\varepsilon}_{k+1}-\boldsymbol{\varepsilon}_{k})(\mathbf{I}-\mathbf{W}),\mathbf{x}_{k+1}-\mathbf{x}^{*}\rangle.
	\end{aligned}
\end{equation}
Since $\gamma<\frac{\lambda_2}{\tau\rho_{A}}$ with $\lambda_2=\lambda_{\min}(\mathbf{W}), \rho_A=\lambda_{\max}(\mathbf{A}^T\mathbf{A})$ $\mathbf{M}=\begin{bmatrix}\mathbf{W}&\mathbf{0}&-\tau\mathbf{A}\\\mathbf{0}&\tau^2\mathbf{L}^{-1}/\psi &\mathbf{0}\\-\tau\mathbf{A}^{T}&\mathbf{0}&\frac{\tau}{\gamma}\mathbf{I}\end{bmatrix} \succ 0$, using property $$\begin{aligned}&2 \left\langle\mathbf{a}-\mathbf{b},\mathbf{a}-\mathbf{c}\right\rangle_{\mathbf{M}}=\left\|\mathbf{a}-\mathbf{c}\right\|_{\mathbf{M}}^{2}-\left\|\mathbf{b}-\mathbf{c}\right\|_{\mathbf{M}}^{2}+\left\|\mathbf{a}-\mathbf{b}\right\|_{\mathbf{M}}^{2}, \end{aligned}$$
from \eqref{tgamma}, we can further deduce that
\begin{equation} \label{midvalue}
	\begin{aligned}
		&\left\|\mathbf{x}_{k+1}-\mathbf{x}^{*}\right\|_{\mathbf{W}}^2+\left\|\mathbf{x}_{k+1}-\mathbf{x}_{k}\right\|_{\mathbf{W}}^2- \left\|\mathbf{x}_{k}-\mathbf{x}^{*}\right\|_{\mathbf{W}}^2\\
		&\quad+\tau^2/\psi \left\|\mathbf{y}_{k+1}-\mathbf{y}^{*}\right\|_{\mathbf{L}^{-1}}^2+\tau^2/\psi \left\|\mathbf{y}_{k+1}-\mathbf{y}_{k}\right\|_{\mathbf{L}^{-1}}^2\\
		&\quad-\tau^2/\psi \left\|\mathbf{y}_{k}-\mathbf{y}^{*}\right\|_{\mathbf{L}^{-1}}^2+\frac{\tau}{\gamma}\left\|\mathbf{z}_{k+1}-\mathbf{z}^{*}\right\|^2\\
		&\quad+\frac{\tau}{\gamma}\left\|\mathbf{z}_{k+1}-\mathbf{z}_{k}\right\|^2+\tau\langle\mathbf{x}_{k}-\mathbf{x}^{*},\mathbf{A}(\mathbf{z}_{k}-\mathbf{z}^{*})\rangle\\
		&\quad-\tau\langle\mathbf{x}_{k+1}-\mathbf{x}^{*},\mathbf{A}(\mathbf{z}_{k+1}-\mathbf{z}^{*})\rangle-\frac{\tau}{\gamma}\left\|\mathbf{z}_{k}-\mathbf{z}^{*}\right\|^2\\
		&\quad-\tau\langle\mathbf{x}_{k+1}-\mathbf{x}_{k},\mathbf{A}(\mathbf{z}_{k+1}-\mathbf{z}_{k})\rangle\\
		&\quad-\tau\langle\mathbf{A}(\mathbf{z}_{k+1}-\mathbf{z}^{*}),\mathbf{x}_{k+1}-\mathbf{x}^{*}\rangle\\
		&\quad+\tau\langle\mathbf{A}(\mathbf{z}_{k}-\mathbf{z}^{*}),\mathbf{x}_{k}-\mathbf{x}^{*}\rangle\\
		&\quad-\tau\langle\mathbf{A}(\mathbf{z}_{k+1}-\mathbf{z}_{k}),\mathbf{x}_{k+1}-\mathbf{x}_{k}\rangle\\
		& = \tau\beta\left\|\Delta\mathbf{z}_{k+1}\right\|^{2}-\frac{2l_{f}L_{f}\tau}{L_{f}+l_{f}}\left\|\mathbf{z}_{k}-\mathbf{z}^{*}\right\|^{2} \\
		&\quad-\left(\frac{2\tau}{L_f+l_f}-\frac\tau\beta\right)\left\|\nabla f(\mathbf{z}_k)-\nabla f(\mathbf{z}^*)\right\|^2\\
		&\quad-2\tau\langle\mathbf{y}_{k+1}-\mathbf{y}^*,\boldsymbol{\varepsilon}_{k+1}\rangle\\
		&\quad+2\psi\langle(\boldsymbol{\varepsilon}_{k+1}-\boldsymbol{\varepsilon}_{k})(\mathbf{I}-\mathbf{W}),\mathbf{x}_{k+1}-\mathbf{x}^{*}\rangle.
	\end{aligned}
\end{equation}
Let $\mathbf{s}_{k}=\left[\mathbf{x}_k^T, \mathbf{y}_k^T, \mathbf{z}_k^T\right]^{T}$ represent the sequence generated by the proposed compressed algorithm \eqref{compact}. The deduction in \eqref{midvalue} can be denoted by the following compact form 
\begin{equation}
	\begin{aligned}
		\|\mathbf{s}_{k+1}-&\mathbf{s}^{*}\|_{\mathbf{M}}^2-\|\mathbf{s}_{k}-\mathbf{s}^{*}\|_{\mathbf{M}}^2+\|\Delta \mathbf{s}_{k+1}\|_{\mathbf{M}}^{2}\\
		& \le \tau\beta\left\|\Delta\mathbf{z}_{k+1}\right\|^{2}-\frac{2l_{f}L_{f}\tau}{L_{f}+l_{f}}\left\|\mathbf{z}_{k}-\mathbf{z}^{*}\right\|^{2} 
		\\
		&-\left(\frac{2\tau}{L_f+l_f}-\frac\tau\beta\right)\left\|\nabla f(\mathbf{z}_k)-\nabla f(\mathbf{z}^*)\right\|^2\\
		&+2\psi\langle(\boldsymbol{\varepsilon}_{k+1}-\boldsymbol{\varepsilon}_{k})(\mathbf{I}-\mathbf{W}),\mathbf{x}_{k+1}-\mathbf{x}^{*}\rangle\\
		&-2\tau\langle\mathbf{y}_{k+1}-\mathbf{y}^*,\boldsymbol{\varepsilon}_{k+1}\rangle,
	\end{aligned}
\end{equation}
which completes the proof.   $\hfill\blacksquare$

\section{Proof of Theorem 1}              

According to \eqref{equation13}, we can obtain the following relation for $\forall \nu >0$: 
\begin{equation}
	\begin{aligned}
	&2\left\|\mathbf{x}_{k+1}-\mathbf{x}^{*}\right\|^{2}-2\left\|\mathbf{x}_{k}-\mathbf{x}^{*}\right\|^{2}\\
	&\le(1+\nu)\|\nabla f(\mathbf{z}_k)-\nabla f(\mathbf{z}^*)\|^2-2\left\|\Delta\mathbf{x}_{k+1}\right\|^2\\
	&+\left(1+\frac{1}{\nu}\right)\frac{1}{\gamma^2}\|\Delta\mathbf{z}_{k+1}\|^2-\left\|\mathbf{x}_k-\mathbf{x}^*\right\|^2.
	\end{aligned}
\end{equation} 
Multiply both sides by $\frac{\epsilon}{2}$, then we have 
\begin{equation} \label{x2}
	\begin{aligned}
		&\epsilon\left\|\mathbf{x}_{k+1}-\mathbf{x}^{*}\right\|^{2}-\epsilon\left\|\mathbf{x}_{k}-\mathbf{x}^{*}\right\|^{2}-\frac{\epsilon}{2}\left\|\mathbf{x}_k-\mathbf{x}^*\right\|^2 \\
		&\le \frac{\epsilon(1+\nu)}{2}\|\nabla f(\mathbf{z}_k)-\nabla f(\mathbf{z}^*)\|^2\\
		&\quad+\left(1+\frac{1}{\nu}\right)\frac{\epsilon}{2\gamma^2}\begin{Vmatrix}\Delta\mathbf{z}_{k+1}\end{Vmatrix}^2 -\epsilon\left\|\Delta\mathbf{x}_{k+1}\right\|^2.
	\end{aligned}
\end{equation}

Moreover, there exist 
$ \tau\mathbf{y}_{k+1}-\tau\mathbf{y}^{*}=\mathbf{W}\Delta \mathbf{x}_{k+1}+\tau(\mathbf{z}_{k}-\mathbf{z}^{*})-\psi (\boldsymbol{\varepsilon}_{k+1}-\boldsymbol{\varepsilon}_{k})(\mathbf{I}-\mathbf{W})
$, using the property $\left\|\mathbf{U}_1+\mathbf{U}_{2}+\mathbf{U}_3\right\|^2\leq \tau'\left\|\mathbf{U}_1\right\|^2+\frac{2\tau'}{\tau'-1}\left[\left\|\mathbf{U}_{2}\right\|^2+\left\|\mathbf{U}_3\right\|^2\right]$ with $\tau'>1$, then it can be deduced that 
$$
\begin{aligned}
	&\tau^{2}\left\|\mathbf{y}_{k+1}-\mathbf{y}^{*}\right\|^{2}-\tau^{2}\left\|\mathbf{y}_{k}-\mathbf{y}^{*}\right\|^{2}\\
	&\leq \tau' \left\|\Delta\mathbf{x}_{k+1}\right\|_{\mathbf{W}^2}^{2}+\tau^{2}\frac{2\tau'}{\tau'-1}\left\|\mathbf{z}_{k}-\mathbf{z}^{*}\right\|^{2}\\
	&\quad +\frac{2\tau'\psi^2}{\tau'-1}\left\|\mathbf{I}-\mathbf{W}\right\|^2\left\|\boldsymbol{\varepsilon}_{k+1}-\boldsymbol{\varepsilon}_{k}\right\|^2-\tau^{2}\left\|\mathbf{y}_{k}-\mathbf{y}^{*}\right\|^{2},
\end{aligned}
$$
multiplying both side with $\epsilon$, then there exists,

\begin{equation} \label{y2}
	\begin{aligned}
		&\epsilon \tau^{2}\left\|\mathbf{y}_{k+1}-\mathbf{y}^{*}\right\|^{2}-\epsilon \tau^{2}\left\|\mathbf{y}_{k}-\mathbf{y}^{*}\right\|^{2}\\
		&\leq \epsilon \tau' \left\|\Delta\mathbf{x}_{k+1}\right\|_{\mathbf{W}^2}^{2}+\epsilon \tau^{2}\frac{2\tau'}{\tau'-1}\left\|\mathbf{z}_{k}-\mathbf{z}^{*}\right\|^{2}\\
		&\quad+\frac{2\epsilon \tau' \psi^2}{\tau'-1}\left\|\mathbf{I}-\mathbf{W}\right\|^2\left\|\boldsymbol{\varepsilon}_{k+1}-\boldsymbol{\varepsilon}_{k}\right\|^2\\
		&\quad-\epsilon \tau^{2}\left\|\mathbf{y}_{k}-\mathbf{y}^{*}\right\|^{2}.
	\end{aligned}
\end{equation}
Substitute (\ref{x2}) and (\ref{y2}) into (\ref{midvalue}), then we have

\begin{equation} \label{finalinequality}
	\begin{aligned}
		&	\|\mathbf{s}_{k+1}-\mathbf{s}^{*}\|_{\mathbf{H}}^{2}-\|\mathbf{s}_{k}-\mathbf{s}^{*}\|_{\mathbf{H}}^{2} \leq-\|\Delta\mathbf{s}_{k+1}\|_{\mathbf{M}-\mathbf{L}_1}^2\\
		&\quad-\|\mathbf{s}_k-\mathbf{s}^*\|_{\mathbf{L}_2}^2 -2\tau\langle\mathbf{y}_{k+1}-\mathbf{y}^*,\boldsymbol{\varepsilon}_{k+1}\rangle\\
		&\quad-\left(\frac{2\tau}{L_f+l_f}-\frac{\tau}{\beta}-\frac{(1+\nu)\epsilon}{2}\right)\|\nabla f(\mathbf{z}_k)-\nabla f(\mathbf{z}^*)\|^2\\
		&\quad+2\psi \langle(\boldsymbol{\varepsilon}_{k+1}-\boldsymbol{\varepsilon}_{k})(\mathbf{I}-\mathbf{W}),\mathbf{x}_{k+1}-\mathbf{x}^{*}\rangle\\
		&\quad+\frac{2\epsilon \tau'\psi^2}{\tau'-1}\left\|\mathbf{I}-\mathbf{W}\right\|^2\left\|\boldsymbol{\varepsilon}_{k+1}-\boldsymbol{\varepsilon}_{k}\right\|^2,
	\end{aligned}
\end{equation}
with  $\mathbf{H}=\begin{bmatrix}\epsilon \mathbf{I}+\mathbf{W}&\mathbf{0}&-\tau\mathbf{A}\\\mathbf{0}&\tau^2(\epsilon \mathbf{I}+\mathbf{L}^{-1}/\psi )&\mathbf{0}\\-\tau\mathbf{A}^{T}&\mathbf{0}&\frac{\tau}{\gamma}\mathbf{I}\end{bmatrix}$, 
\\
$\begin{gathered}\mathbf{L}_{1}=\begin{bmatrix}\epsilon(\tau'\mathbf{W}^2-\mathbf{I})&\mathbf{0}&\mathbf{0}\\\mathbf{0}&\mathbf{0}&\mathbf{0}\\\mathbf{0}&\mathbf{0}&\left(\frac12\left(1+\frac1\nu\right)\frac\tau\gamma+\tau\beta\right)\mathbf{I}\end{bmatrix},\\ 
	\mathbf{L}_{2}=\begin{bmatrix}\frac\epsilon2\mathbf{I}&\mathbf{0}&\mathbf{0}\\0&\epsilon\tau^2\mathbf{I}&\mathbf{0}\\0&0&\left(\frac{2L_fl_f\tau}{L_f+l_f}-\frac{2 \tau'}{\tau'-1}\epsilon\tau^2\right)\mathbf{I}\end{bmatrix}.\end{gathered}$

Recall the compression errors $\boldsymbol{\varepsilon}_{k}=\hat{\mathbf{x}}_{k}-\mathbf{x}_{k}$ and $\boldsymbol{\varepsilon}_{k + 1}=\hat{\mathbf{x}}_{k + 1}-\mathbf{x}_{k + 1}$,  we first address the compression error related terms in equation \eqref{finalinequality}. Scaling these error terms directly leads to  
\begin{equation}   \label{error}
	\begin{aligned}
		&-2\tau\mathbb{E}\left[\langle\mathbf{y}_{k+1}-\mathbf{y}^*,\boldsymbol{\varepsilon}_{k+1}\rangle\right]\\
		&+2\psi\mathbb{E}\left[\langle(\boldsymbol{\varepsilon}_{k+1}-\boldsymbol{\varepsilon}_{k})(\mathbf{I}-\mathbf{W}),\mathbf{x}_{k+1}-\mathbf{x}^{*}\rangle\right]\\
		&+\frac{2\epsilon \tau' \psi^2}{\tau'-1}\mathbb{E}\left[\left\|\mathbf{I}-\mathbf{W}\right\|^2\left\|\boldsymbol{\varepsilon}_{k+1}-\boldsymbol{\varepsilon}_{k}\right\|^2\right] \\
		&=2\psi \mathbb{E}\left[\langle(\mathbf{I}-\mathbf{W})\boldsymbol{\varepsilon}_{k+1},\boldsymbol{\varepsilon}_{k+1}\rangle\right]+2\psi\mathbb{E}\left[\langle(\boldsymbol{\varepsilon}_{k+1}\right.\\
		&\quad\left.-\boldsymbol{\varepsilon}_{k})(\mathbf{I}-\mathbf{W}),(\mathbf{W}-\mathbf{I})\boldsymbol{\varepsilon}_{k}\rangle\right]\\
		&\quad+\frac{2\epsilon \tau' \psi^2}{\tau'-1}\mathbb{E}\left[\left\|\mathbf{I}-\mathbf{W}\right\|^2\left\|\boldsymbol{\varepsilon}_{k+1}-\boldsymbol{\varepsilon}_{k}\right\|^2\right]\\
		 &\le 2\psi\left\|\mathbf{I}-\mathbf{W}\right\|\mathbb{E}\left[\left\|\boldsymbol{\varepsilon}_{k+1}\right\|^2\right]\\
		&\quad+{2\psi\mathbb{E}\left[\langle(\boldsymbol{\varepsilon}_{k+1}-\boldsymbol{\varepsilon}_{k})(\mathbf{I}-\mathbf{W}),(\mathbf{W}-\mathbf{I})\boldsymbol{\varepsilon}_{k}\rangle\right]} \\
		&+\frac{2\epsilon \tau' \psi^2}{\tau'-1}\mathbb{E}\left[\left\|\mathbf{I}-\mathbf{W}\right\|^2\left\|\boldsymbol{\varepsilon}_{k+1}-\boldsymbol{\varepsilon}_{k}\right\|^2\right]\\
	&\le 2\psi \left\|\mathbf{I}-\mathbf{W}\right\| \mathbb{E}\left[\left\|\boldsymbol{\varepsilon}_{k+1}\right\|^2\right]+\psi \left\|{\mathbf{I}-\mathbf{W}}\right\|^2\mathbb{E}\left[\left\|\boldsymbol{\varepsilon}_{k+1}\right\|^2\right] \\
		&\quad+\psi \left\|{\mathbf{I}-\mathbf{W}}\right\|^2\mathbb{E}\left[\left\|\boldsymbol{\varepsilon}_{k}\right\|^2\right]+2\psi \left\|\mathbf{I}-\mathbf{W}\right\|^2
		\mathbb{E}\left[\left\|\boldsymbol{\varepsilon}_{k}\right\|^2\right]\\
		&\quad+\frac{4\epsilon \tau'\psi^2}{\tau'-1}\left\|\mathbf{I}-\mathbf{W}\right\|^2\mathbb{E}\left[\left\|\boldsymbol{\varepsilon}_{k+1}\right\|^2\right]\\
		&\quad+ \frac{4\epsilon \tau'\psi^2}{\tau'-1}\left\|\mathbf{I}-\mathbf{W}\right\|^2\mathbb{E}\left[\left\|\boldsymbol{\varepsilon}_{k}\right\|^2\right]\\
		&\le c_1r_k^2\sigma^2+c_2r_{k+1}^2\sigma^2,
	\end{aligned}
\end{equation}
with $c_1=(3\psi+\frac{4\epsilon \tau'\psi^2}{\tau'-1})\left\|\mathbf{I}-\mathbf{W}\right\|^2, c_2=2\psi\left\|\mathbf{I}-\mathbf{W}\right\|+(\psi+\frac{4\epsilon \tau'\psi^2}{\tau'-1})\left\|\mathbf{I}-\mathbf{W}\right\|^2$. It is obvious that $c_1<c_2$  by coefficient comparison. If the condition  
$\xi<\frac{c_1}{c_2}$ holds, we yields the key contraction relation $c_2\xi^{k+1}<c_1\xi^{k}$.  The first inequality is obtained by using the Cauchy-Schwarz inequality $\langle \mathbf{u}, \mathbf{v} \rangle \le \|\mathbf{u}\|\|\mathbf{v}\|$.  The second inequality is proved by using the  mean inequality relation $\langle \mathbf{a}, \mathbf{b}\rangle \le \|\mathbf{a}\|^2+\|\mathbf{b}\|^2$ and vector squared triangle inequality $\|\mathbf{x}-\mathbf{y}\|^2 \le 2\|\mathbf{x}\|^2 + 2\|\mathbf{y}\|^2$. Let $r_k^2 = h\xi^{k}$, where $\xi\in(0, 1)$ and $h$ is a constant. According to the condition $c_2r_{k + 1}^2 < c_1r_k^2$, it is equivalent to $c_2\xi^{k+1}<c_1\xi^{k}$, we get $1 - \delta<\xi<\frac{c_1}{c_2}$.  In this case, $c_1r_k^2\sigma^2 + c_2r_{k + 1}^2\sigma^2\leq2c_1r_k^2\sigma^2$. 
In  view of (\ref{finalinequality}), it reduces to verifying that $\mathbf{M}>\mathbf{L}_1, \mathbf{L}_{2}>\delta \mathbf{H}, \left(\frac{2\tau}{L_f+l_f}-\frac{\tau}{\beta}-\frac{(1+\nu)\epsilon}{2}\right)\ge 0$. We first prove that  $\mathbf{M}>\mathbf{L}_1$. 
Be definition, it suffices to ensure the matrix condition
$$
\epsilon \mathbf{I}+\mathbf{W}-\epsilon\tau'\mathbf{W}^2>0. 
$$
Since a matrix is positive definite if and only all its eigenvalues are positive, the above condition is equivalent to $\epsilon+\lambda-\epsilon\tau'\lambda^2>0,
	$ 
	for every eigenvalue $\lambda$ of $\mathbf{W}$. Substituting 
	$\mu=1, \epsilon=\gamma \tau, \tau^{\prime}=2, \psi=1$ and letting  $\lambda_2=\lambda_{\min}(\mathbf{W})$, this inequality holds if $\gamma<\frac{\lambda_2}{\tau}$. Recall that $\lambda_{\max }(\mathbf{W})=1$. Then the minimum eigenvalue satisfies  \(\lambda_{\min }(\epsilon \mathbf{I}+\mathbf{W}-\epsilon(1+\) \(\mu) \mathbf{W}^{2})=\lambda_2-\epsilon \mu\). With \(\mu=1\) and \(\nu=2\) , the inequality $ \left(\frac{\tau}{\gamma}-\frac{1}{2}\left(1+\frac{1}{\nu}\right) \frac{\tau}{\gamma}-\tau \beta\right) \mathbf{I} -\tau^{2}\mathbf{A}^{T}\left(\epsilon\mathbf{I}+\mathbf{W}-\epsilon(1+\mu) \mathbf{W}^{2}\right)^{-1}\mathbf{A}>0$ reduces to  $\frac{1}{4\gamma} - \beta > \frac{\tau \rho_{{A}}}{\lambda_2 - \gamma \tau}$. Rearranging terms yields the upper bound  $\gamma<\frac{\lambda_{2}}{4 \lambda_{2}\beta+\tau(4 \rho_{A}+1)}$.  Therefore, if $\gamma$  satisfies $$\gamma <\min\left\{\frac{\lambda_{2}}{\tau}, \frac{\lambda_{2}}{4 \lambda_{2}\beta+\tau(4 \rho_{A}+1)}\right\},$$ we can conclude that  \(\mathbf{M}>\mathbf{L}_{1}\).
	Next, we consider the condition $\left(\frac{2\tau}{L_f+l_f}-\frac{\tau}{\beta}-\frac{(1+\nu)\epsilon}{2}\right)\ge 0$. By substituting  $\mu=1$ and $\nu=2$, and recalling $\epsilon=\gamma\tau$, the above inequality can be ensured by choosing the step-size $\gamma$ such that 
	$$\gamma<\frac{2}{3}\left(\frac{2}{L_{f}+l_{f}}-\frac{1}{\beta}\right).$$ Under this choice, it follows that  $\left(\frac{2\tau}{L_f+l_f}-\frac{\tau}{\beta}-\frac{(1+\nu)\epsilon}{2}\right)>0$. Finally, we proceed to establish the matrix inequality $\mathbf{L}_{2}>\delta \mathbf{H}$ with matrices

$\mathbf{L}_{2}=\begin{bmatrix}\frac\epsilon2\mathbf{I}&\mathbf{0}&\mathbf{0}\\0&\epsilon\tau^2\mathbf{I}&\mathbf{0}\\0&0&\left(\frac{2L_fl_f\tau}{L_f+l_f}-\frac{2 \tau'}{\tau'-1}\epsilon\tau^2\right)\mathbf{I}\end{bmatrix}$,\\
$\delta\mathbf{H}=\begin{bmatrix}\delta(\epsilon \mathbf{I}+\mathbf{W})&\mathbf{0}&-\delta\tau\mathbf{A}\\\mathbf{0}&\delta\tau^2(\epsilon I+\mathbf{L}^{-1}/\psi )&\mathbf{0}\\-\delta\tau\mathbf{A}^{T}&\mathbf{0}&\frac{\tau\delta}{\gamma}\mathbf{I}\end{bmatrix}$, \\
if $\delta < \min\left\{ \frac{\epsilon}{2\epsilon(1+\rho_A) + 2}, \gamma\left( \frac{L_f l_f}{L_f + l_f} - \tau \right),\frac{\epsilon(1-\eta)}{2-\eta} \right\}$ with $\beta=\left(L_{f}+l_{f}\right), \rho_{A}=\lambda_{\max}(\mathbf{A}^{T}\mathbf{A})$, there exists $\mathbf{L}_{2}>\delta \mathbf{H}$.
According to the above derivation results, it follows that: 
\begin{equation}
	\begin{aligned}
		&\mathbb{E}\left[\|\mathbf{s}_{k+1}-\mathbf{s}^*\|_{\mathbf{H}}^2\right]\leq(1-\delta)\mathbb{E}\left[\|\mathbf{s}_k-\mathbf{s}^*\|_{\mathbf{H}}^2\right] + 2c_1r_k^2\sigma^2\\
		&\le (1-\delta)^{k+1}\mathbb{E}\left[\|\mathbf{s}_0-\mathbf{s}^*\|_{\mathbf{H}}^2\right]+\sum_{l=0}^{k}(1-\delta)^{k-l}\xi^{l}2c_1h\sigma^2\\
		&\le (1-\delta)^{k+1}\mathbb{E}\left[\|\mathbf{s}_0-\mathbf{s}^*\|_{\mathbf{H}}^2\right]+c\xi^{k+1},
	\end{aligned}
\end{equation}
where $c=\frac{2c_1h\sigma^2}{\xi-(1-\delta)}$, which completes the proof.  $\hfill\blacksquare$

\section{Proof of Theorem 2}  
Based on the deduction results in \eqref{midvalue}, we can obtain
\begin{equation} \label{contraction2}
	\begin{aligned}
		&-2\tau\mathbb{E}\left[\langle\mathbf{y}_{k+1}-\mathbf{y}^*,\boldsymbol{\varepsilon}_{k+1}\rangle\right]\\
		&+2\psi\mathbb{E}\left[\langle(\boldsymbol{\varepsilon}_{k+1}-\boldsymbol{\varepsilon}_{k})(\mathbf{I}-\mathbf{W}),\mathbf{x}_{k+1}-\mathbf{x}^{*}\rangle\right]\\
		&+\frac{2\epsilon \tau' \psi^2}{\tau'-1}\mathbb{E}\left[\left\|\mathbf{I}-\mathbf{W}\right\|^2\left\|\boldsymbol{\varepsilon}_{k+1}-\boldsymbol{\varepsilon}_{k}\right\|^2\right] \\
		&=2\psi \mathbb{E}\left[\langle(\mathbf{I}-\mathbf{W})\boldsymbol{\varepsilon}_{k+1},\boldsymbol{\varepsilon}_{k+1}\rangle \right]+2\psi\mathbb{E}\left[\langle(\boldsymbol{\varepsilon}_{k+1}\right.\\
		&\quad\left.-\boldsymbol{\varepsilon}_{k})(\mathbf{I}-\mathbf{W}),(\mathbf{W}-\mathbf{I})\boldsymbol{\varepsilon}_{k}\rangle\right]\\
		&\quad+\frac{2\epsilon \tau' \psi^2}{\tau'-1}\left\|\mathbf{I}-\mathbf{W}\right\|^2\mathbb{E}\left[\left\|\boldsymbol{\varepsilon}_{k+1}-\boldsymbol{\varepsilon}_{k}\right\|^2\right] \\
		&\le 2\psi\left\|\mathbf{I}-\mathbf{W}\right\|\mathbb{E}\left[\left\|\boldsymbol{\varepsilon}_{k+1}\right\|^2\right]\\
		&\quad+{2\psi\mathbb{E}\left[\langle(\boldsymbol{\varepsilon}_{k+1}-\boldsymbol{\varepsilon}_{k})(\mathbf{I}-\mathbf{W}),(\mathbf{W}-\mathbf{I})\boldsymbol{\varepsilon}_{k}\rangle\right]} 
		\\
		&\quad+\frac{2\epsilon \tau' \psi^2}{\tau'-1}\left\|\mathbf{I}-\mathbf{W}\right\|^2\mathbb{E}\left[\left\|\boldsymbol{\varepsilon}_{k+1}-\boldsymbol{\varepsilon}_{k}\right\|^2\right] \\
		&\le 2\psi \left\|\mathbf{I}-\mathbf{W}\right\| \mathbb{E}\left[\left\|\boldsymbol{\varepsilon}_{k+1}\right\|^2\right]+\psi \left\|{\mathbf{I}-\mathbf{W}}\right\|^2\mathbb{E}\left[\left\|\boldsymbol{\varepsilon}_{k+1}\right\|^2\right]\\
		&\quad+\psi \left\|{\mathbf{I}-\mathbf{W}}\right\|^2\mathbb{E}\left[\left\|\boldsymbol{\varepsilon}_{k}\right\|^2\right]+2\psi \left\|\mathbf{I}-\mathbf{W}\right\|^2
		\mathbb{E}\left[\left\|\boldsymbol{\varepsilon}_{k}\right\|^2\right]\\
		&\quad+\frac{4\epsilon \tau'\psi^2}{\tau'-1}\left\|\mathbf{I}-\mathbf{W}\right\|^2\left(\mathbb{E}\left[\left\|\boldsymbol{\varepsilon}_{k+1}\right\|^2\right]+\mathbb{E}\left[\left\|\boldsymbol{\varepsilon}_{k}\right\|^2\right]\right)\\
		&\le d_1 C\mathbb{E}\left[\left\|\mathbf{x}_{k+1}-\mathbf{h}_{x}^{k+1}\right\|^2\right]+ d_2 C \mathbb{E}\left[\left\|\mathbf{x}_{k}-\mathbf{h}_{x}^{k}\right\|^2\right],
	\end{aligned}
\end{equation}
with $d_1=d_2=\left(3\psi+\frac{4\epsilon \tau'\psi^2}{\tau'-1}\right)\left\|\mathbf{I}-\mathbf{W}\right\|$ and $\psi, \epsilon \in(0, 1)$.

By exploiting the contractive property of the compression operator $Q_2(x)$, we derive the following relation for the compression error:
\begin{equation}  \label{compressinequality}
	\begin{aligned}
		& \mathbb{E}\left[\left\|\mathbf{x}_{k+1}-\mathbf{h}_{x}^{k+1}\right\| ^{2}\right] \\
		&=\mathbb{E}\left[\left\| \mathbf{x}_{k+1}-\mathbf{x}_{k}+\mathbf{x}_{k}-\mathbf{h}_{x}^{k}-\alpha_{x} Q_2(x)\right\| ^{2}\right] \\ 
		&\leq  \tau_{x}\mathbb{E}\left[\left\| \alpha_{x} (\mathbf{x}_{k}-\mathbf{h}_{x}^{k}-{C}(\mathbf{x}_{k}-\mathbf{h}_{x}^{k})) \right.\right.\\
		&\quad\left.\left.+(1-\alpha_{x})(\mathbf{x}_{k}-\mathbf{h}_{x}^{k})\right\|^{2}\right]+\frac{\tau_{x}}{\tau_{x}-1}\mathbb{E}\left[\left\| \mathbf{x}_{k+1}-\mathbf{x}_{k}\right\| ^{2}\right] \\ 
		& \leq \tau_{x}\alpha_{x} \mathbb{E}\left[\left\|\mathbf{x}_{k}-\mathbf{h}_{x}^{k}-{C}(\mathbf{x}_{k}-\mathbf{h}_{x}^{k})\right\|^2\right]\\
		& \quad+(1-\alpha_{x})\mathbb{E}\left[\left\|\mathbf{x}_{k}-\mathbf{h}_{x}^{k}\right\|^{2}\right]+\frac{\tau_{x}}{\tau_{x}-1}\mathbb{E}\left[\left\| \mathbf{x}_{k+1}-\mathbf{x}_{k}\right\| ^{2}\right]
		\\ 
		& \leq  \tau_{x}\left[\alpha_{x} (1-\delta)+\left(1-\alpha_{x} \right)\right]\mathbb{E}\left[\left\| \mathbf{x}_{k}-\mathbf{h}_{x}^{k}\right\| ^{2}\right] \\
		&\quad +\frac{\tau_{x}}{\tau_{x}-1} \mathbb{E}\left[\left\| \mathbf{x}_{k+1}-\mathbf{x}_{k}\right\| ^{2} \right],
	\end{aligned}
\end{equation}
with the constant $c_{x}= \tau_{x}\left[\alpha_{x} (1-\delta)+\left(1-\alpha_{x} \right)\right]=\tau_{x}(1-\alpha_xr\delta)<1$ where the parameters satisfy  $\frac{\tau_{x}}{\tau_{x}-1}>1$. The inequality stems from the property of compressor $Q_2(\cdot)$. 

Substitute the result in (\ref{contraction2}) into \eqref{finalinequality} and multiply both sides of \eqref{compressinequality} by $a \in (0,1)$, and  add the resulting inequality to \eqref{contraction2}. It follows that

\begin{equation} \label{thm2derivation}
	\begin{aligned}
		&(\epsilon \mathbf{I}+
		\mathbf{W})\left(\mathbb{E}\left[\left\|\mathbf{x}_{k+1}-\mathbf{x}^{*}\right\|^2\right]-\mathbb{E}\left[\left\|\mathbf{x}_{k}-\mathbf{x}^{*}\right\|^2\right]\right)\\
		&+(\epsilon \tau^{2}+ \frac{\tau^2\mathbf{L}^{-1}}{\psi} )\left(\mathbb{E}\left[\left\|\mathbf{y}_{k+1}-\mathbf{y}^{*}\right\|^2\right]-\mathbb{E}\left[\left\|\mathbf{y}_{k}-\mathbf{y}^{*}\right\|^2\right]\right)\\
		&+\frac{\tau}{\gamma}\mathbb{E}\left[\left\|\mathbf{z}_{k+1}-\mathbf{z}^{*}\right\|^2\right]-\frac{\tau}{\gamma}\mathbb{E}\left[\left\|\mathbf{z}_{k}-\mathbf{z}^{*}\right\|^2\right]\\
		&-\tau\mathbb{E}\left[\langle\mathbf{x}_{k+1}-\mathbf{x}^{*},\mathbf{z}_{k+1}-\mathbf{z}^{*}\rangle\right]+\tau\mathbb{E}\left[\langle\mathbf{x}_{k}-\mathbf{x}^{*}, \right.\\
	    &\left.\mathbf{z}_k-\mathbf{z}^*\rangle \right]-\tau\mathbb{E}\left[\langle \mathbf{x}_{k+1}-\mathbf{x}_k, \mathbf{z}_{k+1}-\mathbf{z}_k \rangle \right]\\
		&-\tau\mathbb{E}\left[\langle\mathbf{z}_{k+1}-\mathbf{z}^{*},\mathbf{x}_{k+1}-\mathbf{x}^{*}\rangle\right]+\tau\mathbb{E}\left[\langle\mathbf{z}_{k}-\mathbf{z}^{*},\right.\\
		&\left.\mathbf{x}_{k}-\mathbf{x}^{*}\rangle\right]-\tau\mathbb{E}\left[\langle\mathbf{z}_{k+1}-\mathbf{z}_{k},\mathbf{x}_{k+1}-\mathbf{x}_{k}\rangle\right]\\
		&+a\mathbb{E}\left[\left\|\mathbf{x}_{k+1}-\mathbf{h}_{x}^{k+1}\right\|^2\right]\\
		&\le \left(-\mathbf{W}+\epsilon(\tau' \mathbf{W}^2-\mathbf{I})+\frac{a\tau_x}{\tau_x-1}\right)\mathbb{E}\left[\left\|\Delta\mathbf{x}_{k+1}\right\|^{2}\right]\\
		&\quad-\tau^2{\mathbf{L}^{-1}}/\psi \mathbb{E}\left[\left\|\Delta\mathbf{y}_{k+1}\right\|^2\right]+ac_x \mathbb{E}\left[\left\|\mathbf{x}_{k}-\mathbf{h}_{x}^{k}\right\|^2\right]\\ 
		&\quad+\left(\left(1+\frac{1}{\nu}\right)\frac{\epsilon}{2\gamma^2}+\tau\beta\right)\mathbb{E}\left[\left\|\Delta\mathbf{z}_{k+1}\right\|^{2}\right]\\
		&\quad-\left(\frac{2\tau}{L_f+l_f}-\frac\tau\beta\right)\mathbb{E}\left[\left\|\nabla f(\mathbf{z}_k)-\nabla f(\mathbf{z}^*)\right\|^2\right]\\
		&\quad-\frac{\epsilon(1+\nu)}{2}\mathbb{E}\left[\left\|\nabla f(\mathbf{z}_k)-\nabla f(\mathbf{z}^*)\right\|^2\right]\\
		&\quad-\epsilon\tau^{2}\mathbb{E}\left[\left\|\mathbf{y}_{k}-\mathbf{y}^{*}\right\|^{2}\right]-\frac{\epsilon}{2}\mathbb{E}\left[\left\|\mathbf{x}_k-\mathbf{x}^*\right\|^2\right]\\
		&\quad-\left(\frac{2l_{f}L_{f}\tau}{L_{f}+l_{f}}-\frac{2 \tau'\epsilon\tau^{2}}{\tau'-1}\right)\mathbb{E}\left[\left\|\mathbf{z}_{k}-\mathbf{z}^{*}\right\|^{2}\right]\\
		&\quad+ d_1 C\mathbb{E}\left[\left\|\mathbf{x}_{k+1}-\mathbf{h}_{x}^{k+1}\right\|^2\right]+ d_2 C \mathbb{E}\left[\left\|\mathbf{x}_{k}-\mathbf{h}_{x}^{k}\right\|^2\right].
	\end{aligned}
\end{equation}
Then, based on the above inequality derivations in \eqref{thm2derivation}, we obtain (in a compact form)
\begin{equation} 
	\begin{aligned}
		&	\|\mathbf{s}_{k+1}-\mathbf{s}^{*}\|_{\mathbf{H}}^{2}-\|\mathbf{s}_{k}-\mathbf{s}^{*}\|_{\mathbf{H}}^{2}\\
		&\leq-\|\Delta\mathbf{s}_{k+1}\|_{\mathbf{M}-\mathbf{L}_1}^2-\|\mathbf{s}_k-\mathbf{s}^*\|_{\mathbf{L}_2}^2\nonumber
			\end{aligned}
	\end{equation}
		\begin{equation} \label{finalinequality1}
			\begin{aligned}
		&\quad-\left(\frac{2\tau}{L_f+l_f}-\frac{\tau}{\beta}-\frac{(1+\nu)\epsilon}{2}\right)\|\nabla f(\mathbf{z}_k)-\nabla f(\mathbf{z}^*)\|^2\\
		&\quad+ d_1 C\mathbb{E}\left[\left\|\mathbf{x}_{k+1}-\mathbf{h}_{x}^{k+1}\right\|^2\right]+ d_2 C \mathbb{E}\left[\left\|\mathbf{x}_{k}-\mathbf{h}_{x}^{k}\right\|^2\right].
	\end{aligned}
\end{equation}

Let $\mathbf{M}, \mathbf{H}$ and $\mathbf{L}_2$ be defined as  in Theorem 1. In  light of (\ref{thm2derivation}), it suffices to verify the following three conditions: $$\mathbf{M}>\mathbf{L}_1, \mathbf{L}_{2}>\delta \mathbf{H}, \left(\frac{2\tau}{L_f+l_f}-\frac{\tau}{\beta}-\frac{(1+\nu)\epsilon}{2}\right)\ge 0. $$  Here, the matrix $\mathbf{L}_1$ is revised as \\
 $\begin{gathered}\mathbf{L}_{1}=\begin{bmatrix}\epsilon(\tau'\mathbf{W}^2-\mathbf{I})+\frac{a\tau_x}{\tau_x-1}\mathbf{I} &\mathbf{0}&\mathbf{0}\\\mathbf{0}&\mathbf{0}&\mathbf{0}\\\mathbf{0}&\mathbf{0}&\left(\frac12\left(1+\frac1\nu\right)\frac\tau\gamma+\tau\beta\right)\mathbf{I}\end{bmatrix}\end{gathered}$,\\and selecting  $\frac{a\tau_x}{\tau_x-1}=\frac{\epsilon}{2}$. By choosing the step-size  $\gamma$  such that $\gamma < \frac{2\lambda_2}{3\tau}$, then one immediately has $\epsilon \mathbf{I} + \mathbf{W} - \epsilon\tau^{\prime}\mathbf{W}^2-\frac{\epsilon}{2} \mathbf{I} >0$. According to the Schur complement condition, the matrix inequality $ \left(\frac{\tau}{\gamma}-\frac{1}{2}\left(1+\frac{1}{\nu}\right) \frac{\tau}{\gamma}-\tau \beta\right) \mathbf{I} -\tau^{2}\mathbf{A}^{T}\left(\epsilon\mathbf{I}+\mathbf{W}-\epsilon(1+\mu) \mathbf{W}^{2}-\frac{\epsilon}{2}\mathbf{I}\right)^{-1}\mathbf{A}>0$ holds if there exists $\gamma < \frac{2\lambda_{2}}{8 \lambda_{2}\beta+\tau(8 \rho_{A}+3)}$. Next, substituting  $\mu=1$ and $\nu=2$, the condition $\left(\frac{2\tau}{L_f+l_f}-\frac{\tau}{\beta}-\frac{(1+\nu)\epsilon}{2}\right)>0$ is satisfied by selecting 
 $\gamma<\frac{2}{3}\left(\frac{2}{L_{f}+l_{f}}-\frac{1}{\beta}\right)$ , which further implies $\left(\frac{2\tau}{L_f+l_f}-\frac{\tau}{\beta}-\frac{(1+\nu)\epsilon}{2}\right)>0$. By choosing  $\delta$  the same as in Theorem 1, we have $\mathbf{L}_{2}>\delta \mathbf{H}$. By summarizing the above constraints, the step-size $\gamma$ is required to satisfy
$$
 \gamma < \min\left\{
 \frac{2\lambda_2}{3\tau},\
 \frac{2\lambda_2}{8\lambda_2\beta+\tau(8\rho_A+3)},\
 \frac{2}{3}\left(
 \frac{2}{L_f+l_f} - \frac{1}{\beta}
 \right)
 \right\}
$$
 to ensure $\mathbf{M} > \mathbf{L}_1$, $\mathbf{L}_2 > \delta \mathbf{H}$, and the scalar condition
$\frac{2\tau}{L_f+l_f}-\frac{\tau}{\beta} - \frac{(1+\nu)\epsilon}{2}>0$  hold simultaneously. Consequently, the relation in \eqref{thm2derivation} can be rearranged as follows:

\begin{equation} \label{result3}
	\begin{aligned}
		& \mathbb{E}\left[\|\mathbf{s}_{k+1}-\mathbf{s}^*\|_{\mathbf{H}}^2\right]+a\mathbb{E}\left[\left\|\mathbf{x}_{k+1}-\mathbf{h}_{x}^{k+1}\right\| ^{2}\right]\\
		&\leq(1-\delta)\mathbb{E}\left[\|\mathbf{s}_k-\mathbf{s}^*\|_{\mathbf{H}}^2\right]+a\mathbb{E}\left[\left\|\mathbf{x}_{k+1}-\mathbf{h}_{x}^{k+1}\right\| ^{2}\right] \\
		& \leq(1-\delta)\mathbb{E}\left[\|\mathbf{s}_k-\mathbf{s}^*\|_{\mathbf{H}}^2\right]+ d_1 C\mathbb{E}\left[\left\|\mathbf{x}_{k+1}-\mathbf{h}_{x}^{k+1}\right\|^2\right]\\
		&+ d_2 C \mathbb{E}\left[\left\|\mathbf{x}_{k}-\mathbf{h}_{x}^{k}\right\|^2\right]+ac_x\mathbb{E}\left[\left\| \mathbf{x}_{k}-\mathbf{h}_{x}^{k}\right\| ^{2}\right],\\
	\end{aligned}
\end{equation}
Based on the definition of the Lyapunov function \begin{equation}
	\begin{aligned}
	\mathbb{E}\left[V_{k+1}\right]&=\mathbb{E}\left[\|\mathbf{s}_{k+1}-\mathbf{s}^*\|_{\mathbf{H}}^2\right]\\
	&+(a-d_1C)\mathbb{E}\left[\left\|\mathbf{x}_{k+1}-\mathbf{h}_{x}^{k+1}\right\| ^{2}\right],
	\end{aligned}
\end{equation}
we choose $d_1C<a<\frac{\tau_x-1}{\tau_x}$ and $C<\frac{a(1-c_x)}{d_1+d_2}$ to guarantee the positiveness and convergence of the Lyapunov function, respectively.  By substituting the rearranged result \eqref{result3} into the Lyapunov function definition and simplifying, we first derive the one-step contraction property:
\begin{equation}
	\mathbb{E}\left[V_{k+1}\right]\le (1-\nu) \mathbb{E}\left[V_{k}\right],
\end{equation}
where the contraction coefficient $\nu$ is defined as $1-\nu=\text{max}\{1-\delta, \frac{d_2C+ac_{x}}{a-d_1C}\}$. By iterating the above inequality for $k$ steps, we further obtain the exponential convergence result $\mathbb{E}\left[V_{k+1}\right]\le (1-\nu)^{k}\mathbb{E}\left[V_0\right]$. This completes the convergence proof.  $\hfill\blacksquare$

\section{Proof of Theorem 3}  

Based on the deduction results in \eqref{midvalue}, and using Cauchy-Schwarz inequality $2\langle \mathbf{a},\mathbf{b}\rangle\le \frac{1}{\tau}\left\|\mathbf{a}\right\|^2+\tau\left\|\mathbf{b}\right\|^2$ for any $\tau>0$, we can derive that 
\begin{equation} \label{contraction3}
\begin{aligned}
	&-2\tau\langle\mathbf{y}_{k+1}-\mathbf{y}^*,\boldsymbol{\varepsilon}_{k+1}\rangle\\
	&+2\psi\langle(\boldsymbol{\varepsilon}_{k+1}-\boldsymbol{\varepsilon}_{k})(\mathbf{I}-\mathbf{W}),\mathbf{x}_{k+1}-\mathbf{x}^{*}\rangle\\
	&+\frac{2\epsilon \tau' \psi^2}{\tau'-1}\left\|\mathbf{I}-\mathbf{W}\right\|^2\left\|\boldsymbol{\varepsilon}_{k+1}-\boldsymbol{\varepsilon}_{k}\right\|^2 \\
	&\le \tau^3 \left\|\mathbf{y}_{k+1}-\mathbf{y}^*\right\|^2+\frac{1}{\tau^2}\left\|\boldsymbol{\varepsilon}_{k+1}\right\|^2\\
	&\quad+\frac{2\epsilon \tau' \psi^2}{\tau'-1}\left\|\mathbf{I}-\mathbf{W}\right\|^2\left\|\boldsymbol{\varepsilon}_{k+1}-\boldsymbol{\varepsilon}_{k}\right\|^2 \\	&\quad+2\langle\boldsymbol{\varepsilon}_{k+1}-\boldsymbol{\varepsilon}_{k},\psi(\mathbf{I}-\mathbf{W})(\mathbf{x}_{k+1}-\mathbf{x}^{*})\rangle\\
	& \le \tau^3 \left\|\mathbf{y}_{k+1}-\mathbf{y}^*\right\|^2+\frac{1}{\tau^2}\left\|\boldsymbol{\varepsilon}_{k+1}\right\|^2\\
	&\quad+\frac{2\epsilon \tau' \psi^2}{\tau'-1}\left\|\mathbf{I}-\mathbf{W}\right\|^2\left\|\boldsymbol{\varepsilon}_{k+1}-\boldsymbol{\varepsilon}_{k}\right\|^2 \\
	&\quad+2\langle\boldsymbol{\varepsilon}_{k+1}-\boldsymbol{\varepsilon}_{k},-\tau\Delta \mathbf{y}_{k+1}-\psi(\mathbf{I}-\mathbf{W})\boldsymbol{\varepsilon}_{k+1}\rangle\\
	& \le \tau^3 \left\|\mathbf{y}_{k+1}-\mathbf{y}^*\right\|^2+\frac{1}{\tau^2}\left\|\boldsymbol{\varepsilon}_{k+1}\right\|^2\\
	&\quad+\frac{2\epsilon \tau' \psi^2}{\tau'-1}\left\|\mathbf{I}-\mathbf{W}\right\|^2\left\|\boldsymbol{\varepsilon}_{k+1}-\boldsymbol{\varepsilon}_{k}\right\|^2\\
	&\quad+ 2\langle\boldsymbol{\varepsilon}_{k+1}-\boldsymbol{\varepsilon}_{k},-\tau\Delta \mathbf{y}_{k+1}\rangle\\
	&\quad-2\langle\boldsymbol{\varepsilon}_{k+1}-\boldsymbol{\varepsilon}_{k},\psi(\mathbf{I}-\mathbf{W})\boldsymbol{\varepsilon}_{k+1}\rangle\\
	& \le \tau^3 \left\|\mathbf{y}_{k+1}-\mathbf{y}^*\right\|^2+\frac{1}{\tau^2}\left\|\boldsymbol{\varepsilon}_{k+1}\right\|^2\\
	&\quad+\frac{2\epsilon \tau' \psi^2}{\tau'-1}\left\|\mathbf{I}-\mathbf{W}\right\|^2\left\|\boldsymbol{\varepsilon}_{k+1}-\boldsymbol{\varepsilon}_{k}\right\|^2 \\
	&\quad+ \frac{1}{\tau^2}\left\|\boldsymbol{\varepsilon}_{k+1}-\boldsymbol{\varepsilon}_{k}\right\|^2+\tau^3\left\|\Delta \mathbf{y}_{k+1}\right\|^2\\
	&\quad-2\langle\boldsymbol{\varepsilon}_{k+1}-\boldsymbol{\varepsilon}_{k},\psi(\mathbf{I}-\mathbf{W})\boldsymbol{\varepsilon}_{k+1}\rangle.
	\end{aligned}
\end{equation}

According to the property $ 	\tau^3\left\|\mathbf{y}_{k+1}-\mathbf{y}^{*}\right\|^2+\tau^3\left\|\Delta \mathbf{y}_{k+1}\right\|^2=\tau^3\left\|\mathbf{y}_{k+1}-\mathbf{y}_{k}+\mathbf{y}_{k}-\mathbf{y}^{*}\right\|^2+\tau^3\left\|\Delta \mathbf{y}_{k+1}\right\|^2$ together with the inequality $\|\mathbf{a}+\mathbf{b}\|^2\le 2\|\mathbf{a}\|^2+2\|\mathbf{b}\|^2$, one has $\tau^3\left\|\mathbf{y}_{k+1}-\mathbf{y}^{*}\right\|^2+\tau^3\left\|\Delta \mathbf{y}_{k+1}\right\|^2 \le 3\tau^3\left\|\Delta \mathbf{y}_{k+1}\right\|^2+ 2\tau^3\\\left\|
\mathbf{y}_{k}-\mathbf{y}^{*}\right\|^2$. Substituting the above relation into \eqref{contraction3}, we derive the following  inequality:
\begin{equation} \label{compressor3ineuqlity}
	\begin{aligned}
		&-2\tau\langle\mathbf{y}_{k+1}-\mathbf{y}^*,\boldsymbol{\varepsilon}_{k+1}\rangle\\
		&+2\psi\langle(\boldsymbol{\varepsilon}_{k+1}-\boldsymbol{\varepsilon}_{k})(\mathbf{I}-\mathbf{W}),\mathbf{x}_{k+1}-\mathbf{x}^{*}\rangle\\
		&+\frac{2\epsilon \tau' \psi^2}{\tau'-1}\left\|\mathbf{I}-\mathbf{W}\right\|^2\left\|\boldsymbol{\varepsilon}_{k+1}-\boldsymbol{\varepsilon}_{k}\right\|^2 \\
		&\le 3\tau^3\left\|\Delta \mathbf{y}_{k+1}\right\|^2+ 2\tau^3 \left\|\mathbf{y}_{k}-\mathbf{y}^{*}\right\|^2 \frac{1}{\tau^2}\left\|\boldsymbol{\varepsilon}_{k+1}\right\|^2\\
		&\quad+\left(\frac{1}{\tau^2}+\frac{2\epsilon \tau' \psi^2}{\tau'-1}\left\|\mathbf{I}-\mathbf{W}\right\|^2\right)\left\|\boldsymbol{\varepsilon}_{k+1}-\boldsymbol{\varepsilon}_{k}\right\|^2\\
		&\quad-2\langle\boldsymbol{\varepsilon}_{k+1}-\boldsymbol{\varepsilon}_{k},\psi(I-\mathbf{W})\boldsymbol{\varepsilon}_{k+1}\rangle\\
		&\le 3\tau^3\left\|\Delta \mathbf{y}_{k+1}\right\|^2+ 2\tau^3 \left\|\mathbf{y}_{k}-\mathbf{y}^{*}\right\|^2 \\
		&\quad+ \frac{1}{\tau^2}\left\|\boldsymbol{\varepsilon}_{k+1}\right\|^2+ \frac{2}{\tau^2}\left\|\boldsymbol{\varepsilon}_{k+1}\right\|^2+\frac{2}{\tau^2}\left\|\boldsymbol{\varepsilon}_{k} \right\|^2\\
		&\quad+\frac{4\epsilon \tau' \psi^2}{\tau'-1}\left\|\mathbf{I}-\mathbf{W}\right\|^2\left(\left\|\boldsymbol{\varepsilon}_{k+1}\right\|^2+\left\|\boldsymbol{\varepsilon}_{k}\right\|^2\right) \\
		&\quad+ 2\psi\|\boldsymbol{\varepsilon}_{k + 1}\|^2+2\psi\|\boldsymbol{\varepsilon}_{k}\|^2+\psi \|\mathbf{I}-\mathbf{W}\|^2\|\boldsymbol{\varepsilon}_{k+1}\|^2.
	\end{aligned}
\end{equation} 

To facilitate subsequent convergence analysis, we define two positive constant coefficients: $\theta_1=\frac{2}{\tau^2}+\frac{4\epsilon \tau' \psi^2}{\tau'-1}\|\mathbf{I}-\mathbf{W}\|^2+2\psi, \theta_2=\frac{3}{\tau^2}+\frac{4\epsilon \tau' \psi^2}{\tau'-1}\|\mathbf{I}-\mathbf{W}\|^2+2\psi+ \psi \|\mathbf{I}-\mathbf{W}\|^2$. It is obvious that $\theta_1<\theta_2$  by coefficient comparison.
Let $r_k^2 = h\xi^{k}$ with contraction factor $\xi\in(1-\nu, 1)$. Then the condition  
 $\xi<\frac{\theta_2}{\theta_1}$ always holds, which yields the key contraction relation $\theta_1\xi^{k+1}<\theta_2\xi^{k}$. Combining with the derived result in \eqref{compressor3ineuqlity}, we have
\begin{equation}
	\begin{aligned}
		&\|\mathbf{s}_{k+1}-\mathbf{s}^{*}\|_{\mathbf{H}}^{2}-\|\mathbf{s}_{k}-\mathbf{s}^{*}\|_{\mathbf{H}}^{2} \\ &\leq-\|\Delta\mathbf{s}_{k+1}\|_{\mathbf{M}-\mathbf{\Phi}_1}^2-\|\mathbf{s}_k-\mathbf{s}^*\|_{\mathbf{\Phi}_2}^2 \\
		&\quad+\theta_1 m_{k}^2\sigma^2+ \theta_2m_{k+1}^2\sigma^2-\left(\frac{2\tau}{L_f+l_f}-\frac{\tau}{\beta}\right.\\
		&\left.\quad-\frac{(1+\nu){\epsilon}}{2}\right){\|\nabla f(\mathbf{z}_k)-\nabla f(\mathbf{z}^*)\|^2},
	\end{aligned}
\end{equation}
where the matrices are defined as \\ 
$\begin{gathered}\mathbf{\Phi}_{1}=\begin{bmatrix}\epsilon(\tau'\mathbf{W}^2-\mathbf{I})&\mathbf{0}&\mathbf{0}\\\mathbf{0}&{3\tau^3}&\mathbf{0}\\\mathbf{0}&\mathbf{0}&\left(\frac12\left(1+\frac1\nu\right)\frac\tau\gamma+\tau\beta\right)\mathbf{I}\end{bmatrix}, \\
	\mathbf{\Phi}_{2}=\begin{bmatrix}\frac\epsilon2\mathbf{I}&\mathbf{0}&\mathbf{0}\\0&(\epsilon\tau^2-{2\tau^3})\mathbf{I}&\mathbf{0}\\0&0&\left(\frac{2L_fl_f\tau}{L_f+l_f}-\frac{2 \tau'}{\tau'-1}\epsilon\tau^2\right)\mathbf{I}\end{bmatrix}.\end{gathered}\\$\\

To ensure the negative definiteness of the contraction terms and establish a stable recursive decay, it suffices to verify the following three conditions:
$ \mathbf{M}>\mathbf{\Phi}_1, \mathbf{\Phi}_2>\upsilon \mathbf{H}$ and $\left(\frac{2\tau}{L_f+l_f}-\frac{\tau}{\beta}-\frac{(1+\nu)\epsilon}{2}\right)>0$. These matrix inequalities can be  transformed into tractable constraints on the algorithm parameters as follows: (i) From the second diagonal block of $ \mathbf{M}>\mathbf{\Phi}_1$, we require
$\tau^2\mathbf{L}^{-1}/{\psi}-3\tau^3\mathbf{I}>0$, which yields $\psi<\frac{1}{3\tau(1-\eta)}$. Together with the step-size condition $\gamma <\min\left\{\frac{\lambda_{2}}{\tau}, \frac{\lambda_{2}}{4 \lambda_{2}\beta+\tau(4 \rho_{A}+1)}, \frac{2}{3}\left(\frac{1}{L_{f}+l_{f}}\right) \right\}$, the inequalities  $\mathbf{M}>\mathbf{\Phi}_1$ and $\left(\frac{2\tau}{L_f+l_f}-\frac{\tau}{\beta}-\frac{(1+\nu)\epsilon}{2}\right)>0$ are satisfied; (ii) From $\mathbf{\Phi}_2>\upsilon \mathbf{H}$, we derive the matrix inequality $\epsilon\tau^2-{2\tau^3}-\upsilon (\tau^2\epsilon \mathbf{I}+\tau^2\mathbf{L}^{-1}/\psi)>0$, which gives the feasible range  $0 < \upsilon < \frac{(\epsilon - 2\tau)}{\epsilon+3\tau}$ with $\gamma>2$. Combing with $\upsilon < \min\left\{ \frac{\epsilon}{2\epsilon(1+\rho_A) + 2}, \gamma\left( \frac{L_f l_f}{L_f + l_f} - \tau \right),\frac{\epsilon(1-\eta)}{2-\eta} \right\}$, the inequality $\mathbf{\Phi}_2>\upsilon \mathbf{H}$ holds. Under the above parameter choices, the recursive relation reduces to the standard linear contraction form:
\begin{equation}
	\begin{aligned}
		&\|\mathbf{s}_{k+1}-\mathbf{s}^{*}\|_{\mathbf{H}}^{2}\le(1-\upsilon) \|\mathbf{s}_{k}-\mathbf{s}^{*}\|_{\mathbf{H}}^{2}+2\theta_2r_k^2\sigma^2 \\
		&\le (1-\upsilon)^{k+1}\|\mathbf{s}_0-\mathbf{s}^*\|_{\mathbf{H}}^2+\sum_{l=0}^{k}(1-\upsilon)^{k-l}\xi^{l}2\theta_2h\sigma^2\\
		&\le (1-\upsilon)^{k+1}\|\mathbf{s}_0-\mathbf{s}^*\|_{\mathbf{H}}^2+\varpi\xi^{k+1},
	\end{aligned}
\end{equation}
where $\varpi=\frac{2\theta_2h\sigma^2}{\xi-(1-\upsilon)}$, which completes the proof.   $\hfill\blacksquare$

 \end{document}